\documentclass[a4paper,11pt]{amsart}

\usepackage[english]{babel}
\usepackage[latin1]{inputenc}
\usepackage{amsfonts,amsthm,amssymb,amsmath,latexsym,wasysym,stmaryrd,mathrsfs}
\usepackage{multirow,rotating,subfigure,picins,graphicx,bigstrut,lscape,slashbox}
\usepackage{hyperref}
\usepackage[centerlast]{caption}
\usepackage[usenames]{color}
\usepackage[all,dvips,color]{xy}
\usepackage{fancyhdr}
\usepackage{indentfirst}

\graphicspath{{Figures/}}

\newtheorem{theo}{Theorem}[section]
\newtheorem{lemme}[theo]{Lemma}
\newtheorem{prop}[theo]{Proposition}
\newtheorem{cor}[theo]{Corollary}
\newtheorem{defi}[theo]{Definition}

\newtheorem{question}[theo]{Question}

\newtheorem{remarque}[theo]{Remark}

\def\Q{{\mathbb Q}}
\def\R{{\mathbb R}}
\def\Z{{\mathbb Z}}

\def\C{{\mathcal C}}

\def\H{{\mathcal H}}

\def\T{{\mathcal T}}

\def\ie{{\it i.e. }}

\def\fract#1/#2{\hbox{\leavevmode
  \kern.1em \raise .25ex \hbox{\the\scriptfont0 $#1$}\kern-.1em }\big/
  {\hbox{\kern-.15em \lower .5ex \hbox{\the\scriptfont0 $#2$}} }}

\newcommand{\dessin}[2]{
  \vcenter{\hbox{\includegraphics[height=#1]{#2}}}}
\newcommand{\dessinH}[2]{
  \vcenter{\hbox{\includegraphics[width=#1]{#2}}}}
\newcommand{\func}[3]{
  #1\colon #2\longrightarrow#3}

\title{Khovanov homology and star--like isotopies}
\author{Benjamin \textsc{Audoux}}
\address{Laboratoire E.Picard\\
UFR MIG\\
Universit\'e Paul Sabatier\\
118 route de Narbonne\\
31062 Toulouse\\
France}
\email{audoux@picard.ups-tlse.fr}

\begin{document}

\maketitle

\begin{abstract}
A star--like isotopy for oriented links in $3$--space is an isotopy which uses only Reidemeister moves which correspond to the following singularities of planar curves~:
$$
\dessin{.6cm}{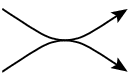}\ \ , \ \ \dessin{.9cm}{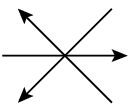}\ \ , \ \ \dessin{.9cm}{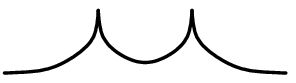}\ \ , \ \ \dessin{.9cm}{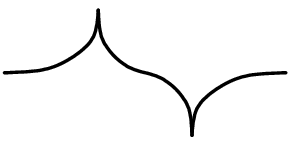}.
$$
We define a link polynomial derived from the Jones polynomial which is, in general, only invariant under star--like isotopies and we categorify it.
\end{abstract}

\section{Introduction}

It is now well known that many polynomial invariants for links and graphs admit a categorification.
We can mention, for instance, the Jones polynomial (\cite{Khovanov}), the Alexander polynomial (\cite{Ozsvath},\cite{Szabo}), the HOMFLYPT polynomial (\cite{Rozansky}), the dichromatic polynomial (\cite{Helme}) or the refinement of the Jones polynomial for braid--like isotopies (\cite{Audoux2}).
Mostly, these invariants are defined for link diagrams and then invariance under Reidemeister moves is proven.
However, in \cite{Rozansky} and \cite{Szabo}, the invariance under Reidemeister IIb moves, \ie Reidemeister II moves corresponding to the following singularity
$$
\dessin{.6cm}{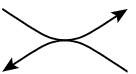},
$$
with opposite tangent vectors, cannot be achieved.
Nevertheless, they can be avoided by using only braid presentations of links. 
On the other hand, the categorification defined in \cite{Rozansky} can be seen as a generalization of the homology given by M. Khovanov in \cite{Khovanov} which is well defined for all diagrams.\\

In \cite{Audoux2}, we have defined a bigraded refinement $Kh_{br}$ of Khovanov invariant $Kh$ which shed light on internal mechanism of invariance for the latter.
As a matter of fact, braid-like Khovanov homology is only invariant under Reidemeister moves which correspond locally to braid isotopies.
In particular, it excludes IIb and Markov moves.
In the construction, the bigrading can be collapsed to a single one.
Via a spectral sequence, the braid-like homology converges then to the usual Khovanov one.\\

In this paper, we give a second refinement of Jones polynomials and its categorification $Kh_{st}$ for another restricted notion of isotopies of link diagrams, called star-like isotopies.
Markov moves are still not allowed, but pairs of Reidemeister I moves on a strand are.
Star-like II moves are the same as braid-like ones but star-like III moves are precisely Reidemeister III moves which are not braid-like.\\
Unfortunately, to date, there is no known geometrical interpretation for star-like isotopies but it gives rise to an homology which is bigraded and which collapses also to Khovanov homology via a spectral sequence (see Fig. \ref{fig:Calcul} for an example).\\

\begin{figure}[t]
$$
\xymatrix @!0 @C=3.5cm @R=3.3cm {
& *+[F-:<5pt>]{\dessin{1.3cm}{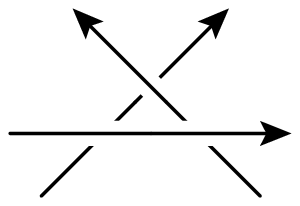}} \ar|{\dessin{.35cm}{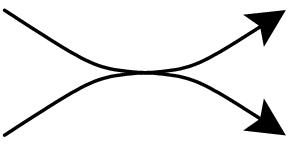}}@{<->}[r] \ar|{\dessin{.35cm}{sing1g}}@{<->}[dr] &
*+[F-:<5pt>]{\dessin{1.3cm}{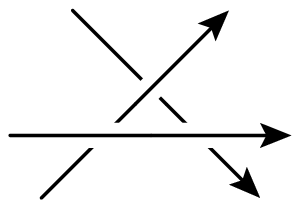}} \ar|{\dessin{.35cm}{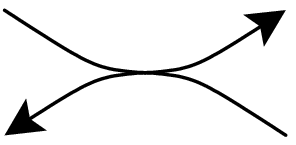}}@{<->}[dr] & \\
*+[F-:<5pt>]{\dessin{1.3cm}{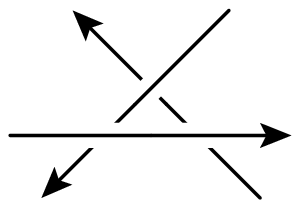}} \ar|{\dessin{.35cm}{sing1bisg}}@{<->}[ur] \ar|{\dessin{.35cm}{sing1bisg}}@{<->}[r] \ar|{\dessin{.35cm}{sing1bisg}}@{<->}[dr] &
*+[F-:<5pt>]{\dessin{1.3cm}{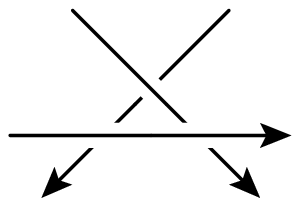}} \ar|{\dessin{.35cm}{sing1g}}@{<->}[ur] \ar|{\dessin{.35cm}{sing1g}}@{<->}[dr] &
*+[F-:<5pt>]{\dessin{1.3cm}{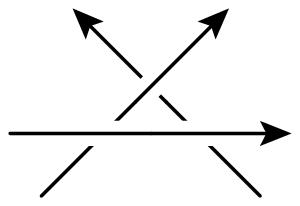}} \ar|{\dessin{.35cm}{sing1bisg}}@{<->}[r] &
*+[F-:<5pt>]{\dessin{1.3cm}{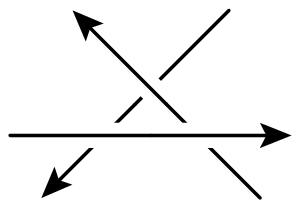}} \\
& *+[F-:<5pt>]{\dessin{1.3cm}{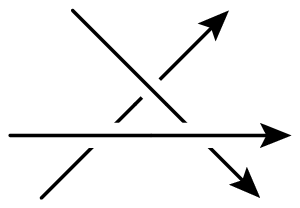}} \ar|{\dessin{.35cm}{sing1g}}@{<->}[ur] \ar|{\dessin{.35cm}{sing1g}}@{<->}[r] &
*+[F-:<5pt>]{\dessin{1.3cm}{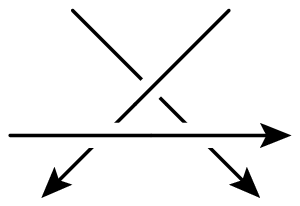}} \ar|{\dessin{.35cm}{sing1bisg}}@{<->}[ur]& \\
}
$$
\caption{Reidemeister III moves graph~: {\footnotesize Each vertex is a Reidemeister move of type III, edges are labeled by Reidemeister moves of type II. An arrow between two vertices means that they can be replaced one by the other using the labeling Reidemeister II move. The mirror operation corresponds to reversing horizontally the graph.}}
\label{diag}
\end{figure}

Besides the enlightening of mechanisms of invariance in Khovanov homology, star-like isotopies are also attractive to study relations between Reidemeister moves.
It is proven in chapter 1.2 of \cite{Fiedler2} that, if an application defined on oriented link diagrams is invariant under the Reidemeister moves corresponding to $\dessin{.36cm}{sing1}$ and under one of the six braid--like Reidemeister moves of type III, then it is also invariant under all the other braid--like Reidemeister III moves.
As shown in Figure \ref{diag}, it is not anymore the case for star--like isotopies.
It points out a dissymmetry in Khovanov construction since a direct proof of invariance can be achieved only for one of the two star-like Reidemeister III moves.
This move depends whether we consider Khovanov homology or cohomology.
Invariance under the other type of move can nonetheless be obtained using a duality argument.\\

\begin{defi}\label{premdef}
Two oriented links are called \emph{star--like} isotopic if their diagrams are related by a sequence of Reidemeister moves of type II and III which correspond to the following singularities of planar curves~:
$$
\dessin{.6cm}{sing1} \hspace{2cm} \dessin{.9cm}{sing2}.
$$
Moreover we allow arbitrary pairs of Reidemeister I moves on the same arc of a diagram.\\
The moves allowed in star--like isotopies are called \emph{star--like Reidemeister moves}.
\end{defi}

The name ``star--like'' comes from the form of the corresponding singularity in star-like Reidemeister III moves.\\

Using so-called surfaces with pulleys, a picture can be given for the construction of star-like Khovanov homology.
It also gives a picture for the braid-like case, as well as the case of links in $I$-bundle (\cite{Asaeda}).

\begin{figure}[h]
\hspace{-2.2cm}
\begin{minipage}{1.0\linewidth}
  \fbox{\begin{xy} (-31,63)*{\textrm{
          $Kh_{st}\left(\dessin{1.5cm}{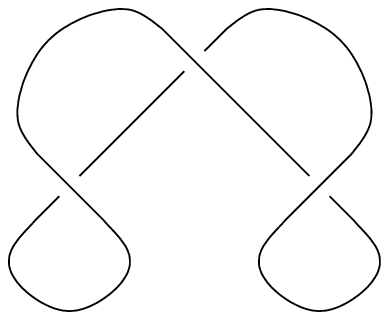}\right) = $
          \begin{tabular}{c|c|c}
            \backslashbox{\hspace{.3cm}$j$}{\hspace{-.5cm}$i$} & 1 & 0 \\
            \hline
            3{\huge \strut} & $\Z_{(0)}$ & \\
            \hline
            1{\huge \strut} & & $\Z_{(2)}\oplus\Z_{(0)}\oplus\Z_{(-2)}$
          \end{tabular}
        }}="Star"; (42,25)*{\textrm{
          \begin{tabular}{c|c|c}
            \backslashbox{\hspace{-.1cm}$j+k$}{\hspace{-.5cm}$i$} & 1 & 0 \\
            \hline
            3{\Large \strut} & $\Z$ & $\Z$ \\
            \hline
            1{\Large \strut} & & $\Z$ \\
            \hline
            -1{\Large \strut} & & $\Z$ 
          \end{tabular}
        }}="CollStar"; (-48,0)*{\textrm{
          $Kh\left(\dessin{1.5cm}{K}\right) = $
          \begin{tabular}{c|c}
            \backslashbox{\hspace{.3cm}$j$}{\hspace{-.5cm}$i$} & 0 \\
            \hline
            1{\Large \strut} & $\Z$ \\
            \hline
            -1{\Large \strut} & $\Z$
          \end{tabular}
        }}="Khovanov"; (47,-20)*{\textrm{
          \begin{tabular}{c|c|c|c}
            \backslashbox{\hspace{-.1cm}$j+k$}{\hspace{-.5cm}$i$} & 1 & 0 & -1 \\
            \hline
            5{\Large \strut} & $\Z$ & $\Z$ & \\
            \hline
            3{\Large \strut} & $\Z^2$ & $\Z^2$ & \\
            \hline
            1{\Large \strut} & $\Z$ & $\Z^3$ & $\Z$ \\
            \hline
            -1{\Large \strut} & & $\Z^2$ & $\Z$ \\
            \hline
            -3{\Large \strut} & & $\Z$ & $\Z$
          \end{tabular}
        }}="CollBraid"; (2,-75)*{\textrm{
          $Kh_{br}\left(\dessin{1.5cm}{K}\right) = $
          \begin{tabular}{c|c|c|c}
            \backslashbox{\hspace{.3cm}$j$}{\hspace{-.5cm}$i$} & 1 & 0 & -1 \\
            \hline
            3{\huge \strut} & $\Z_{(2)}\oplus\Z^2_{(0)}\oplus\Z_{(-2)}$ & &\\
            \hline
            \multirow{2}{*}{1} & & {\huge \strut}$\Z_{(4)}\oplus\Z^2_{(2)}\oplus\Z^3_{(0)}$ &\\
            && {\huge \strut}$\oplus\Z^2_{(-2)}\oplus\Z_{(-4)}$ &\\
            \hline
            -1{\huge \strut} & & & $\Z_{(2)}\oplus\Z_{(0)}\oplus\Z_{(-2)}$ \\
          \end{tabular}
        }}="Braid";

      \ar@/_1.4pc/"Star"!D!<2.4cm,-.3cm>;"CollStar"!UL!<-.3cm,-1cm>_(.35){\textrm{grading
          collapse}}
      \ar@{~>}"CollStar"!DL!<-.3cm,.5cm>;"Khovanov"!R!<.3cm,.5cm>^{\
        \textrm{sequence}}_{\textrm{spectral}\ }
      \ar@{~>}"CollBraid"!UL!<-.3cm,-1.5cm>;"Khovanov"!R!<.3cm,-.2cm>^{\textrm{sequence}\
      }_{\ \textrm{spectral}\ }
      \ar@/^1.4pc/"Braid"!U!<-.9cm,.3cm>;"CollBraid"!DL!<-.3cm,1cm>^(.35){\textrm{grading
          collapse}}
    \end{xy}
  }
\end{minipage}

  \caption{A computation for restricted Khovanov homologies: {\scriptsize for each bigraded copy of $\Z$, the $k$--grading is given between parenthesis.}}
  \label{fig:Calcul}
\end{figure}

\begin{figure}[h]
$$
\begin{array}{ccccc}
\dessin{.85cm}{sing3} & \ \ e.g. \ \ & \dessin{.85cm}{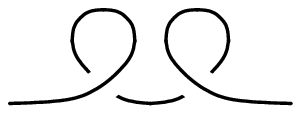} & \longleftrightarrow & \dessin{.85cm}{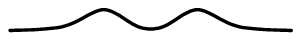}\\
& & & \vdots & \\
& & \dessin{.85cm}{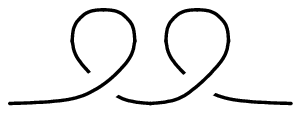} & \longleftrightarrow & \dessin{.85cm}{Istr1}\\
& & \dessin{.85cm}{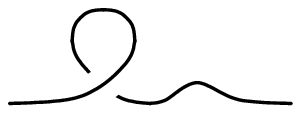} & \longleftrightarrow & \dessin{.85cm}{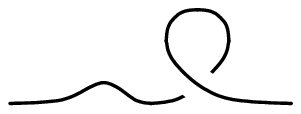}\\
& & & \vdots & \\[5mm]
\dessin{.85cm}{sing4} & \ \ e.g. \ \ & \dessin{.85cm}{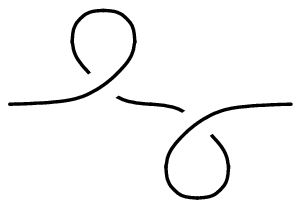} & \longleftrightarrow & \dessin{.85cm}{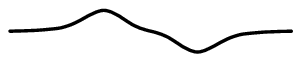}\\
& & & \vdots & \\
& & \dessin{.85cm}{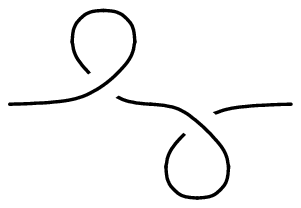} & \longleftrightarrow & \dessin{.85cm}{Istr2}\\
& & \dessin{.85cm}{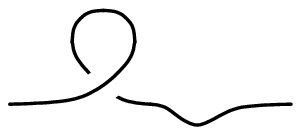} & \longleftrightarrow & \dessin{.85cm}{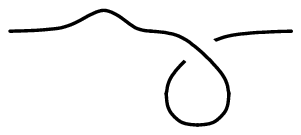}\\
& & & \vdots & \\
\end{array}
$$
\caption{Pairs of Reidemeister I moves allowed in star--like isotopies} \label{Ipair}
\end{figure}

\section{Main results}

\begin{figure}[h]
$$
\xymatrix@!0 @R=.8cm @C=3.5cm {
& \dessin{.85cm}{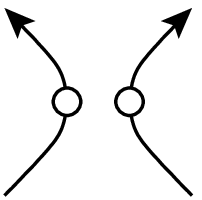} \ \ \textrm{Seifert points}\\
\dessin{.85cm}{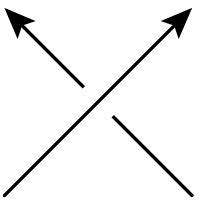} \ar[ur]!L \ar[dr]!L & \\
& \dessin{.85cm}{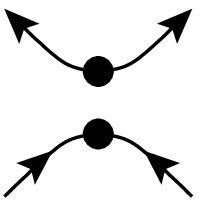} \ \ \textrm{break points}\\
}
$$

\label{fig:Seifert}
\caption{Break and Seifert points}
\end{figure}

As in the braid--like case, our starting point is the Kauffman bracket for oriented links in $3$--space. We use Kauffman's notations and terminology (see \cite{Kauffman}).\\
Each Kauffman state $s$ of an oriented diagram $D$ has a piecewise smooth structure which is induced by the orientation of $D$. The crosssings of $D$ give rise to special points of $s$. If in such a point, the piecewise orientation changes, we say that it is a \emph{break point} (compare with \cite{Audoux2}). The number of break points on any circle of $s$ is always even. The remaining special points are called \emph{Seifert points}. In Seifert points, the orientations of the two pieces fit together (see Fig. \ref{fig:Seifert}).\\

\begin{defi}
A circle in a Kauffman state is of type $d$ if $\big(\frac{1}{2} \#(\textrm{break points}) + \# (\textrm{Seifert points})\big)$ is odd and of type $h$ otherwise.
\end{defi}

Let $C$ denote the set of all configurations, possibly empty, of unoriented circles embedded in the plane up to isotopy.

\begin{defi}
For any Kauffman state $s$, we define
\begin{eqnarray*}
\sigma(s) &  = & \# (A\textrm{--}smoothings) - \# (A^{-1}\textrm{--}smoothings),\\
d(s) & = & \#(d\textrm{--}circles),\\
h(s) & = & \#(h\textrm{--}circles),\\
\end{eqnarray*}
and $c(s) \in C$ which is the configuration of only the $h$--circles.
\end{defi}

Let $\Gamma$ denote the $\Z[A,A^{-1}]$--module over $C$.

\begin{defi}
For any diagram $D$, the bracket $\langle D \rangle_{st} \in \Gamma$ is defined by
$$
\langle D \rangle_{st}=\sum_{\substack{s\textrm{ Kauffman}\\[.1cm] \textrm{state of }D}} A^{\sigma(s)}(-A^2 - A^{-2})^{d(s)} c(s).
$$
The polynomial $V_{st}(D)$ is then defined by
$$
V_{st}(D) = (-A)^{-3w(D)}\langle D \rangle_{st},
$$
where $w(D)$ denotes the writhe of $D$.
\end{defi}

From this definition, it follows immediatly~:
\begin{prop}\label{premprop}
Let $D$ be a diagram and $v$ a crossing of $D$. Let $D_0$ (resp. $D_1$) be the diagram obtained from $D$ by smoothing $v$ in the $A$--fashion (resp. $A^{-1}$--fashion) and where break and Seifert points are added in the natural way. Then, we have
$$
\langle D \rangle_{st} = A \langle D_0 \rangle_{st} + A^{-1} \langle D_1 \rangle_{st}.
$$
Moreover, we have
$$
\langle D \ \dessin{.4cm}{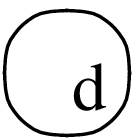} \rangle_{st} = (-A^2 - A^{-2}) \langle D \rangle_{st}.
$$
\end{prop}

\begin{theo}\label{theo1}
The polynomial $V_{st}$ is invariant under star--like isotopies.
\end{theo}

\begin{remarque}
If we denote Seifert points by circles, one easily sees that the polynomial $V_{st}$ satisfies the following skein relations~:
$$
A^4 V_{st}\left( \dessin{.45cm}{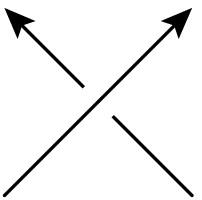} \right) - A^{-4} V_{st} \left( \dessin{.45cm}{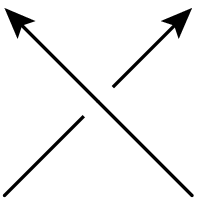} \right) = (A^2 - A^{-2})V_{st}\left( \dessin{.45cm}{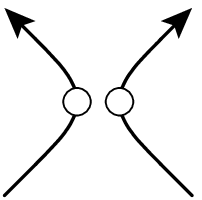} \right).
$$
But $V_{st}\left( \dessin{.45cm}{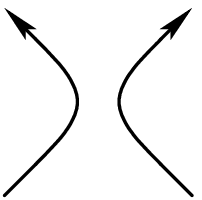} \right) \neq V_{st}\left( \dessin{.45cm}{skeinliss} \right)$ and, consequently, $V_{st}$ does not satisfy Jones skein relation.
\end{remarque}

\begin{remarque}
Replacing $c(s)$ by $X^{h(s)}$ defines a map
$$
\func{\chi}{C}{\Z [A,A^{-1},X]}.
$$
The invariant $V_{st}$ contains certainly more information than $\chi(\langle \ . \ \rangle_{st})$, but it turns out that the latter can be categorified contrary to the former.
\end{remarque}

\begin{theo}\label{theo2}
For any diagram $D$, there are graded chain complexes defined on the same graded module $(\C_{i,j,k}(D))_{i,j,k\in \Z}$ but with two different differentials $d$ and $d'$ of, respectively, tridegree $(-1,0,0)$ and $(1,0,0)$ such that their homology groups, respectively denoted $\H_{***}(D)$ and $\H'_{***}(D)$, are invariant under star--like isotopy.
Moreover, the polynomial $\chi (V_{st})$ is given as the bigraded Euler characteristics of these homologies \ie
\begin{eqnarray*}
\chi\big(V_{st}(D)\big) & = & \sum_{i,j,k} (-1)^i(-A^2)^jX^k dim_\Q \big(\H_{i,j,k}(D)\otimes_\Z \Q\big)\\
& = & \sum_{i,j,k} (-1)^i(-A^2)^jX^k dim_\Q \big(\H'_{i,j,k}(D)\otimes_\Z \Q\big)
\end{eqnarray*}
\end{theo}

\begin{prop}\label{coho}
For any diagram $D$ and any $i,j,k\in \Z$, the groups $\H_{i,j,k}(D)$, $\H'_{-i,-j,-k}(\overline{D})$ and $\H^{-i,-j,-k}(\overline{D})$ are isomorphic, where $\overline{D}$ denotes the miror image of $D$ and $\H^{***}$ the cohomology groups associated to $d$.
\end{prop}

As noticed at the end of section $9$ in \cite{Asaeda}, it follows thus from the Universal Coefficient theorem~:

\begin{cor}\label{mirroir}
For any diagram $D$ and any $i,j,k \in \Z\begin{scriptsize}\begin{footnotesize}\begin{small}\end{small}\end{footnotesize}\end{scriptsize}$, the groups $\H_{i,j,k}(D)$ and $\big(\H_{-i,-j,-k}(\overline{D})/\T_{-i,-j,-k}(\overline{D})\big)\oplus \T_{-i+1,-j,-k}(\overline{D})$ are isomorphic where $\T_{***}(D)$  denote the torsion part of $\H_{***}(D)$.
\end{cor}

\begin{remarque}
Of course, the last two propositions still hold when swapping the role of $d$ and $d'$.
\end{remarque}

\begin{question}
Is it possible to categorify the original polynomial $V_{st}$ ?
\end{question}

\section{Proof of theorem \ref{theo1}}

From now on, we denote break points by dots and Seifert points by small circles on the diagram. Brackets are not always written.\\

\begin{lemme}\label{lem1}
The polynomial $V_{st}$ is invariant under Reidemeister moves of type II with equal tangent direction.
\end{lemme}

\begin{proof}
When calculating the Kauffman bracket, we obtain
$$
\dessin{.75cm}{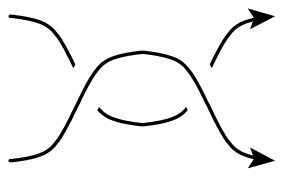} \ = \ A^2 \dessin{.75cm}{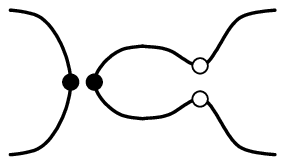} \ + \ \dessin{.75cm}{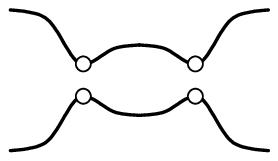} \ + \ \dessin{.75cm}{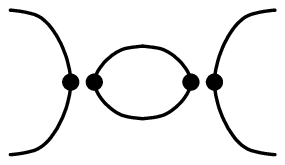}\ + \  A^{-2} \dessin{.75cm}{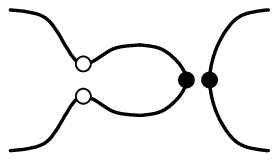} .
$$
The circle is of type $d$ and the usual identification $d = -A^2 - A^{-2}$ implies invariance. 
\end{proof}

\begin{remarque}
Like in the braid--like case, we can not achieve invariance under the other type of Reidemeister II moves. Actually, the calculation gives
$$
\dessin{.75cm}{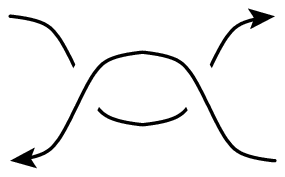} \ = \ A^2 \dessin{.75cm}{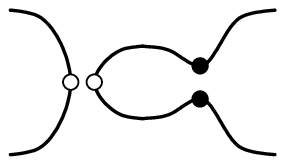} \ + \ \dessin{.75cm}{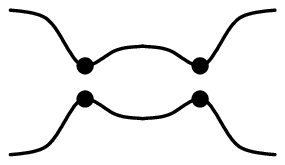}  \ + \ \dessin{.75cm}{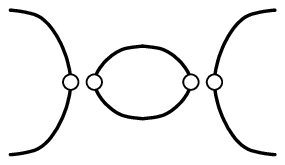}\ + \ A^{-2} \dessin{.75cm}{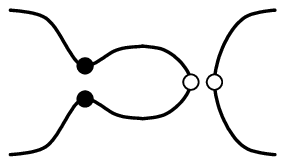},
$$
and the diagrams are completely different with respect to the types $d$ or $h$ of the circles.
\end{remarque}

\begin{lemme}\label{lem2}
The polynomial $V_{st}$ is invariant under star--like Reidemeister moves of type III.
\end{lemme}

\begin{figure}[h]
$$
\begin{array}{ccc}
\dessin{1.1cm}{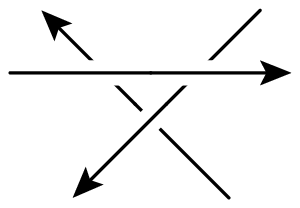} & \hspace{.7cm} \longrightarrow  \hspace{.7cm} & \dessin{1.1cm}{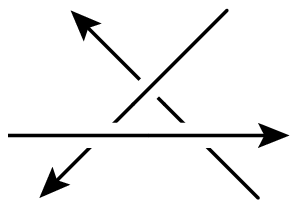}\\
& \textrm{move } III_g &\\[5mm]
\dessin{1.1cm}{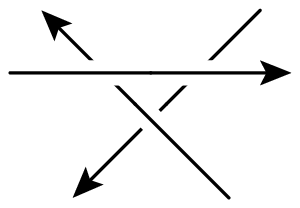} & \hspace{.7cm} \longrightarrow  \hspace{.7cm} & \dessin{1.1cm}{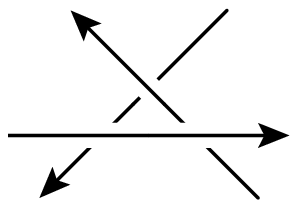}\\
& \textrm{move }III_h &\\
\end{array}
$$
\caption{The two star--like Reidemeister moves of type III }\label{2cas}
\end{figure}

\begin{proof}
We will prove the invariance of $V_{st}$ by the two star--like Reidemeister moves of type III, shown in Figure \ref{2cas}.\\
The left hand side and the right hand side of the first one give respectively the following eight contributions to the bracket~:
$$
\hspace{-.9cm}
(1)
\left\{\textrm{
    \begin{tabular}{p{1.9cm}p{.6cm}p{1.9cm}}
      \begin{tabular}{p{.3cm}p{1.3cm}}
        \centering{$A$} & $\dessin{1.05cm}{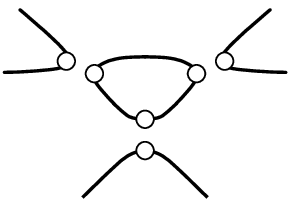}$ 
      \end{tabular}
      &&
      \begin{tabular}{p{.3cm}p{1.3cm}}
        \centering{$A^3$} & $\dessin{1.05cm}{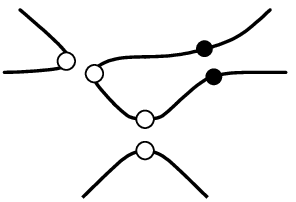}$ 
      \end{tabular}
      \\[.6cm]
      \begin{tabular}{p{.3cm}p{1.3cm}}
        \centering{$A^{-1}$} & $\dessin{1.05cm}{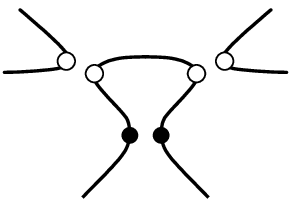}$ 
      \end{tabular}
      &&
      \begin{tabular}{p{.3cm}p{1.3cm}}
        \centering{$A^{-1}$} & $\dessin{1.05cm}{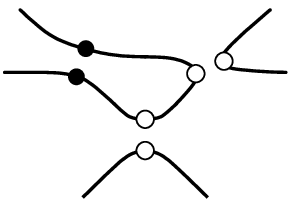}$ 
      \end{tabular}
    \end{tabular}}
\right.
\hspace{.9cm}
(1')
\left\{\textrm{
    \begin{tabular}{p{5.2cm}}
      \centering{\begin{tabular}{p{.3cm}p{1.3cm}}
          \centering{$A^{-1}$} & $\dessin{1.05cm}{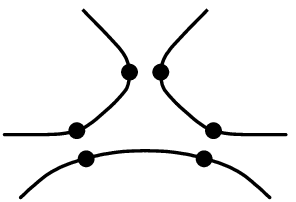}$ 
        \end{tabular}}
    \end{tabular}}   
\right.
$$
$$
\hspace{-.9cm}
(2)
\left\{\textrm{
    \begin{tabular}{p{5.2cm}}
      \centering{\begin{tabular}{p{.3cm}p{1.3cm}}
          \centering{$A^{-1}$} & $\dessin{1.05cm}{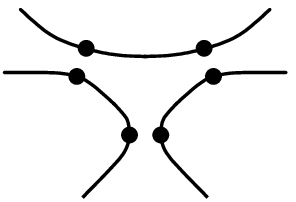}$ 
        \end{tabular}}
    \end{tabular}}   
\right.
\hspace{1cm}
(2')
\left\{\textrm{
    \begin{tabular}{p{1.9cm}p{.6cm}p{1.9cm}}
      \begin{tabular}{p{.3cm}p{1.3cm}}
        \centering{$A$} & $\dessin{1.05cm}{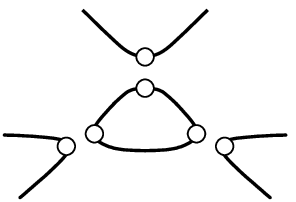}$ 
      \end{tabular}
      &&
      \begin{tabular}{p{.3cm}p{1.3cm}}
        \centering{$A^3$} & $\dessin{1.05cm}{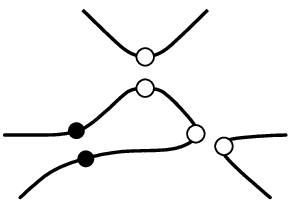}$ 
      \end{tabular}
      \\[.6cm]
      \begin{tabular}{p{.3cm}p{1.3cm}}
        \centering{$A^{-1}$} & $\dessin{1.05cm}{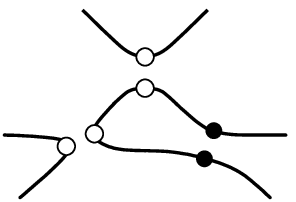}$ 
      \end{tabular}
      &&
      \begin{tabular}{p{.3cm}p{1.3cm}}
        \centering{$A^{-1}$} & $\dessin{1.05cm}{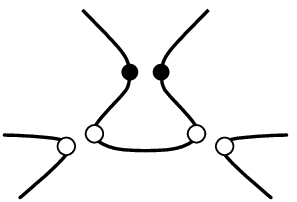}$ 
      \end{tabular}
    \end{tabular}}
\right.
$$
\vspace{.3cm}
$$
\hspace{-.9cm}
(3)
\left\{\textrm{
\begin{tabular}{p{1.9cm}p{.6cm}p{1.9cm}}
  \begin{tabular}{p{.3cm}p{1.3cm}}
    \centering{$A$} & \centering{$\dessin{1.05cm}{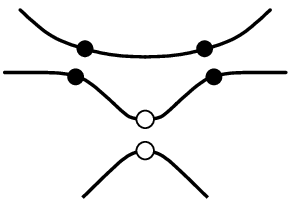}$} 
  \end{tabular}
  &&
  \begin{tabular}{p{.3cm}p{1.3cm}}
    \centering{$A^{-3}$} & \centering{$\dessin{1.05cm}{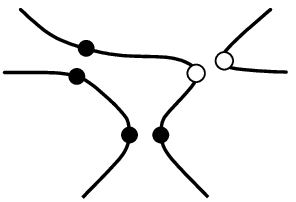}$} 
  \end{tabular}
  \\[.6cm]
  \multicolumn{3}{c}{\begin{tabular}{p{.3cm}p{1.3cm}}
      \centering{$A$} & $\dessin{1.05cm}{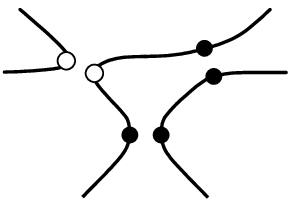}$ 
    \end{tabular}}
\end{tabular}}
\right.
\hspace{1cm}
(3')
\left\{\textrm{
\begin{tabular}{p{1.9cm}p{.6cm}p{1.9cm}}
  \begin{tabular}{p{.3cm}p{1.3cm}}
    \centering{$A$} & \centering{$\dessin{1.05cm}{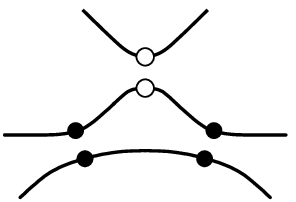}$} 
  \end{tabular}
  &&
  \begin{tabular}{p{.3cm}p{1.3cm}}
    \centering{$A^{-3}$} & \centering{$\dessin{1.05cm}{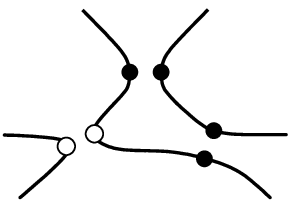}$} 
  \end{tabular}
  \\[.6cm]
  \multicolumn{3}{c}{\begin{tabular}{p{.3cm}p{1.3cm}}
      \centering{$A$} & $\dessin{1.05cm}{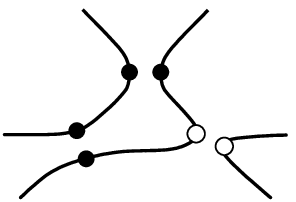}$ 
    \end{tabular}}
\end{tabular}}
\right.
$$

The closed connected components arising in $(1)$ and $(2')$ are of type $d$. This implies that $(1)=(1')$, $(2)=(2')$ and $(3)=(3')$ when summing their elements to evaluate the star-like Kauffman bracket. Moreover, the configurations of $h$--circles in $\R^2$ are the same.\\

The left hand side and the right hand side of the second move give respectively the following contributions~:
$$
\hspace{-.9cm}
(1)
\left\{\textrm{
    \begin{tabular}{p{1.9cm}p{.6cm}p{1.9cm}}
      \begin{tabular}{p{.3cm}p{1.3cm}}
        \centering{$A^{-1}$} & $\dessin{1.05cm}{III11}$ 
      \end{tabular}
      &&
      \begin{tabular}{p{.3cm}p{1.3cm}}
        \centering{$A$} & $\dessin{1.05cm}{III12}$ 
      \end{tabular}
      \\[.6cm]
      \begin{tabular}{p{.3cm}p{1.3cm}}
        \centering{$A^{-1}$} & $\dessin{1.05cm}{III13}$ 
      \end{tabular}
      &&
      \begin{tabular}{p{.3cm}p{1.3cm}}
        \centering{$A$} & $\dessin{1.05cm}{III14}$ 
      \end{tabular}
    \end{tabular}}
\right.
\hspace{.9cm}
(1')
\left\{\textrm{
    \begin{tabular}{p{5.2cm}}
      \centering{\begin{tabular}{p{.3cm}p{1.3cm}}
          \centering{$A$} & $\dessin{1.05cm}{III1s1}$ 
        \end{tabular}}
    \end{tabular}}   
\right.
$$
$$
\hspace{-.9cm}
(2)
\left\{\textrm{
    \begin{tabular}{p{5.2cm}}
      \centering{\begin{tabular}{p{.3cm}p{1.3cm}}
          \centering{$A$} & $\dessin{1.05cm}{III21}$ 
        \end{tabular}}
    \end{tabular}}   
\right.
\hspace{1cm}
(2')
\left\{\textrm{
    \begin{tabular}{p{1.9cm}p{.6cm}p{1.9cm}}
      \begin{tabular}{p{.3cm}p{1.3cm}}
        \centering{$A^{-1}$} & $\dessin{1.05cm}{III2s1}$ 
      \end{tabular}
      &&
      \begin{tabular}{p{.3cm}p{1.3cm}}
        \centering{$A$} & $\dessin{1.05cm}{III2s2}$ 
      \end{tabular}
      \\[.6cm]
      \begin{tabular}{p{.3cm}p{1.3cm}}
        \centering{$A$} & $\dessin{1.05cm}{III2s3}$ 
      \end{tabular}
      &&
      \begin{tabular}{p{.3cm}p{1.3cm}}
        \centering{$A^{-3}$} & $\dessin{1.05cm}{III2s4}$ 
      \end{tabular}
    \end{tabular}}
\right.
$$
\vspace{.3cm}
$$
\hspace{-.9cm}
(3)
\left\{\textrm{
\begin{tabular}{p{1.9cm}p{.6cm}p{1.9cm}}
  \begin{tabular}{p{.3cm}p{1.3cm}}
    \centering{$A^3$} & \centering{$\dessin{1.05cm}{III31}$} 
  \end{tabular}
  &&
  \begin{tabular}{p{.3cm}p{1.3cm}}
    \centering{$A^{-1}$} & \centering{$\dessin{1.05cm}{III32}$} 
  \end{tabular}
  \\[.6cm]
  \multicolumn{3}{c}{\begin{tabular}{p{.3cm}p{1.3cm}}
      \centering{$A^{-1}$} & $\dessin{1.05cm}{III33}$ 
    \end{tabular}}
\end{tabular}}
\right.
\hspace{1cm}
(3')
\left\{\textrm{
\begin{tabular}{p{1.9cm}p{.6cm}p{1.9cm}}
  \begin{tabular}{p{.3cm}p{1.3cm}}
    \centering{$A^3$} & \centering{$\dessin{1.05cm}{III3s1}$} 
  \end{tabular}
  &&
  \begin{tabular}{p{.3cm}p{1.3cm}}
    \centering{$A^{-1}$} & \centering{$\dessin{1.05cm}{III3s2}$} 
  \end{tabular}
  \\[.6cm]
  \multicolumn{3}{c}{\begin{tabular}{p{.3cm}p{1.3cm}}
      \centering{$A^{-1}$} & $\dessin{1.05cm}{III3s3}$ 
    \end{tabular}}
\end{tabular}}
\right.
$$
We conclude like in the first case.
\end{proof}

\begin{lemme}\label{lem3}
The polynomial $V_{st}$ is invariant under star--like Reidemeister I moves.
\end{lemme}

\begin{figure}[h]
$$
\begin{array}{ccc}
\dessin{1.1cm}{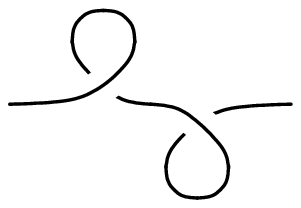} & \hspace{.7cm} \longrightarrow  \hspace{.7cm} & \dessin{1.1cm}{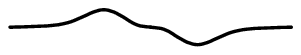}\\
\end{array}
$$
\caption{First star--like Reidemeister move of type I} \label{1stImove}
\end{figure}

\begin{figure}[h]
$$
\begin{array}{ccc}
\dessin{1.1cm}{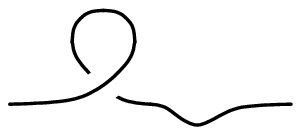} & \hspace{.7cm} \longrightarrow  \hspace{.7cm} & \dessin{1.1cm}{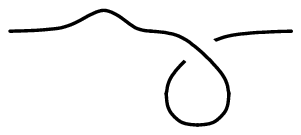}\\
\end{array}
$$
\caption{Second star--like Reidemeister move of type I} \label{2ndImove}
\end{figure}

\begin{proof}
We will only consider two cases. All other cases are analoguous and are left to the reader.

The left hand side of Figure \ref{1stImove} gives the following contributions~:
\begin{eqnarray*}
\dessin{.9cm}{Ipair6} & = & A^2 \dessin{.9cm}{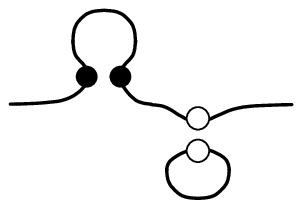} + \dessin{.9cm}{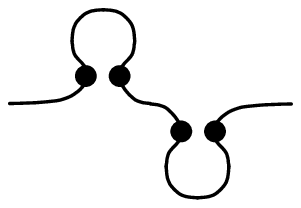} + \dessin{.9cm}{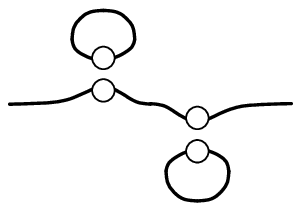} + A^{-2} \dessin{.9cm}{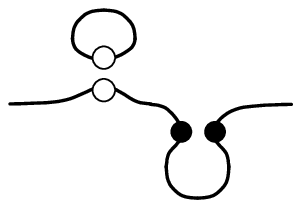}\\
& = & (A^2d + 1 + d^2 + A^{-2}d) \dessin{.9cm}{Istr2}\\
& = & \dessin{.9cm}{Istr2}.
\end{eqnarray*}

Similarly, the left hand side of Figure \ref{2ndImove} gives the contribution~:
\begin{eqnarray*}
-A^3\left( A \dessin{.9cm}{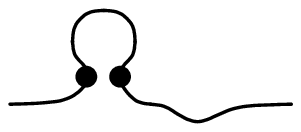} + A^{-1}\dessin{.9cm}{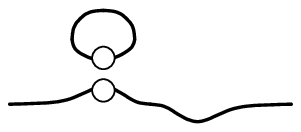}\right) & = & -A^3\left( (A + A^{-1}d) \dessin{.9cm}{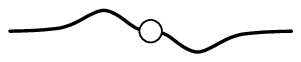}\right)\\
& = & \dessin{.9cm}{I25},
\end{eqnarray*}
whereas its right hand side gives~:
\begin{eqnarray*}
-A^{-3}\left( A \dessin{.9cm}{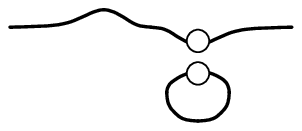} + A^{-1}\dessin{.9cm}{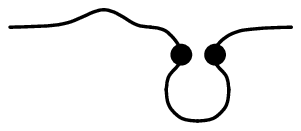}\right) & = & -A^{-3}\left( (Ad + A^{-1}) \dessin{.9cm}{I25}\right)\\
& = & \dessin{.9cm}{I25}.
\end{eqnarray*}

\end{proof}

Theorem \ref{theo1} follows directly from Lemmata \ref{lem1}, \ref{lem2} and \ref{lem3}.

\section{Proof of theorem \ref{theo2}}

Let $D$ be a diagram with $n$ crossings.
\subsection{Definition of homologies}

\begin{defi}
A enhanced state $S$ of $D$ is a Kauffman state $s$ of $D$ enhanced by an assignment of a plus or a minus sign to each of the circles of $s$.
\end{defi}

Setting $X=-H^2 - H^{-2}$ and according to \cite{Viro}, we rewrite our lightened bracket as
\begin{eqnarray}\label{bracket}
\chi\big(\langle D \rangle_{st}\big) = \sum_{\substack{S\textrm{ enhanced}\\[.1cm] \textrm{state of }D}} (-1)^{\frac{\sigma(S) - w(K)}{2}}(-A^2)^{\frac{\sigma(S)-3w(K)}{2}+\tau_d(S)}(-H^2)^{\tau_h(S)},
\end{eqnarray}
where $\tau_d$ (resp. $\tau_h$) is the difference between the number of plus and minus assigned to the $d$--circles (resp. $h$--circles).

\begin{defi}
For any enhanced state $S$ of $D$, we define
$$
\begin{array}{l}
i(S) = \frac{\sigma(S) - w(K)}{2}\\[2mm]
j(S) = \frac{\sigma(S)-3w(K)}{2}+\tau_d(S)\\[2mm]
k(S) = \tau_h(S).
\end{array}
$$ 
Then, for $i,j,k \in \Z$, we can define $\C_{i,j,k}(D)$ to be the $\Z$--module spanned by
$$
\{S\textrm{ enhanced state of }D\ |\ i(S) = i, j(S)=j, k(S)=k\}.
$$
\end{defi}

\begin{defi}
Let $S$ and $S'$ be two enhanced states of $D$ and $v$ a crossing of $D$. We define an incidence number $[S:S']_v$ by

\begin{itemize}
\item[-] $[S:S']_v=1$ if the following four conditions are satisfied~:
\begin{enumerate}
\item[i)] $v$ is $A$--smoothed in $S$ but $A^{-1}$--smoothed in $S'$;
\item[ii)] all the other crossings are smoothed the same way in $S$ and $S'$;
\item[iii)] common circles of $S$ and $S'$ are labeled the same way;
\item[iv)] $j(S) = j(S')$ and $k(S) = k(S')$;
\end{enumerate}
\item[-] otherwise, $[S:S']_v=0$.
\end{itemize}
\end{defi}

Now, we assign an ordering of the crossings of $D$.

\begin{defi}
We define a differential of tridegree $(-1,0,0)$ on $\C_{***}(D)$ by
$$
d(S) =  \sum_{v\textrm{ crossing in }D} d_v(S)
$$
where the partial differential $d_v$ is defined by
$$
d_v(S) = (-1)^{t^-_{v,S}} \sum_{S'\textrm{ state of }D} [S,S']_v S',
$$
with $t^-_{v,S}$ the number of $A^{-1}$--smoothed crossings in S labeled with numbers greater than the label of $v$.
\end{defi}

\begin{defi}
In the same way, we define a differential of tridegree $(1,0,0)$ on $\C_{***}(D)$ by
$$
d'(S) =  \sum_{v\textrm{ crossing in }D} d'_v(S)
$$
where the partial differential $d'_v$ is defined by
$$
d'_v(S) = (-1)^{t^+_{v,S}} \sum_{S'\textrm{ state of }D} [S',S]_v S',
$$
with $t^+_{v,S}$ the number of $A$--smoothed crossings in S labeled with numbers greater than the label of $v$.
\end{defi}

The differential $d(S)$ (resp. $d'(S)$) is the alternating sum of enhanced states obtained by switching one $A$--smoothed (resp. $A^{-1}$--smoothed) crossing of $S$ into a $A^{-1}$--smoothing (resp. $A$--smoothing) and by locally relabeling in such a way that $j$ and $k$ are preserved. Moreover the merge of two $d$--circles or two $h$--circles always gives a $d$--circle and the merge of a $d$--circle with an $h$--circle always an $h$--circle. Thus, one can explicit, for instance, the action of $d_v$, as done in Figure 8.

\footnotetext[1]{depending on the type of the resulting circles}
\begin{lemme}
Every two distinct partial differentials $d_v$ and $d_{v'}$(resp. $d'_v$ and $d'_{v'}$) anti--commute. Thus $d^2=0$ (resp. $d'^2=0$) and $d$ (resp. $d'$) is a differential. 
\end{lemme}
\begin{proof}
One only has to check cases from $1A$--$1D$ to $5A$--$5D$ of the proof of Lemma $4.4$ in \cite{Asaeda}. 
\end{proof}

\begin{lemme}
The homologies $\H(D)$ and $\H'(D)$ do not depend on the ordering of the crossings.
\end{lemme}

\begin{figure}[t]
$$
\begin{array}{ccc}
\begin{array}{ccc}
\dessin{1cm}{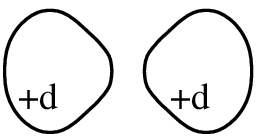}& \longrightarrow &zero\\[7mm]
\dessin{1cm}{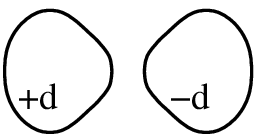}& \longrightarrow  &\dessin{1cm}{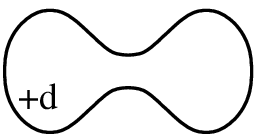}\\[7mm]
\dessin{1cm}{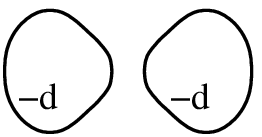}& \longrightarrow  &\dessin{1cm}{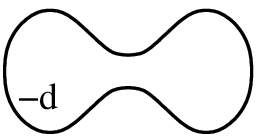}\\[7mm]
\dessin{1cm}{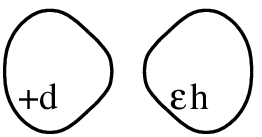}&  \longrightarrow  &zero\\[7mm]
\end{array}& \ \ &
\begin{array}{ccc}
\dessin{1cm}{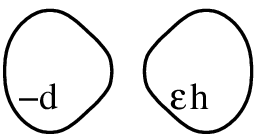}& \longrightarrow  &\dessin{1cm}{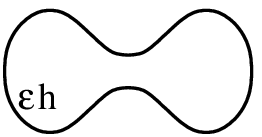}\\[7mm]
\dessin{1cm}{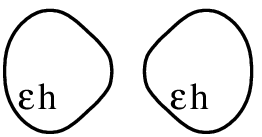}& \longrightarrow &zero\\[7mm]
\dessin{1cm}{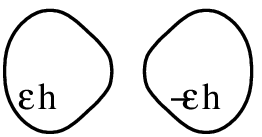}& \longrightarrow  &\dessin{1cm}{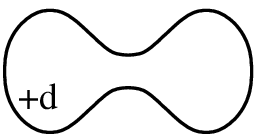}\\[7mm]
\end{array}
\end{array}
$$
$$
\begin{array}{ccl}
\dessin{1cm}{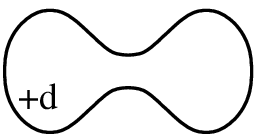}& \longrightarrow  &
\left\{\begin{array}{cl}
\dessin{1cm}{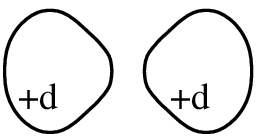}&\footnotemark[1]\\[7mm]
zero&
\end{array}\right.
\\[12mm]
\dessin{1cm}{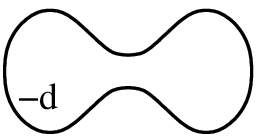}& \longrightarrow &
\left\{ \begin{array}{cl}
\dessin{1cm}{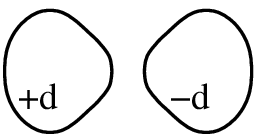}+\dessin{1cm}{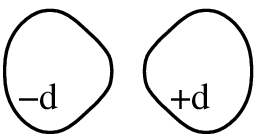}& \footnotemark[1]\\[7mm]
\dessin{1cm}{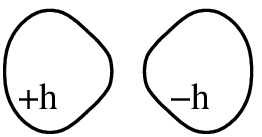}+\dessin{1cm}{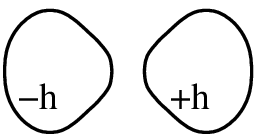}&\\
\end{array} \right.
\\[16mm]
\dessin{1cm}{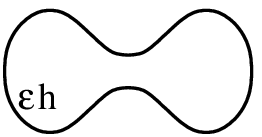}& \longrightarrow  &\dessin{1cm}{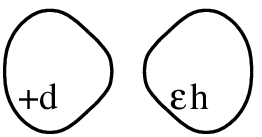}\\
\end{array}
$$
\label{fig:Diff}
\caption{Definition of the partial differential $d_v$}
\end{figure}

\begin{proof}
This follows from the proof of theorem 6.1 in \cite{Asaeda}.
\end{proof}

\begin{prop}
The Euler characteristics of $\H(D)$ and $\H'(D)$ give the same polynomial $V_{st}(D)$.
\end{prop}
\begin{proof}
This is clear from equation \ref{bracket} and the definition of the chain complexes.
\end{proof}

\subsection{Proof of Proposition \ref{coho}}

For any chain complex $(C_{***},\delta)$ where $\delta$ is of tridegree $(a,0,0)$ with $a \in \Z$, we can define its dual complex $(C^{***},\delta^*)$ by
$$
\begin{array}{l}
\forall i,j,k \in \Z,C^{i,j,k}=Hom(C_{i,j,k},\Z);\\[3mm]
\forall \phi \in C^{***},\delta^*(\phi)= \phi \circ \delta.
\end{array}
$$
The map $\delta^*$ is then of tridegree $(-a,0,0)$.

\begin{lemme}
The map
$$
\psi : \vcenter{\hbox{\xymatrix{
\cdots \ar[r]^(.4){d'} & \C'_{i,j,k}(D) \ar[r]^(.45){d'} \ar[d]^\psi & \C'_{i+1,j,k}(D) \ar[r]^(.63){d'} \ar[d]^\psi& \cdots\\
\cdots \ar[r]^(.4){d^*} & \C^{i,j,k}(D) \ar[r]^(.46){d^*} & \C^{i+1,j,k}(D) \ar[r]^(.6){d^*} & \cdots,\\
}}}
$$
defined on any enhanced state $S$ by $\psi(S)=(-1)^{u(S)}S^*$ where $S^*$ is the dual element of $S$ and $u$ is defined by
$$
u(S)=\sum_{\substack{v_i\textrm{ crossing}\\[.1cm]A\textrm{--smoothed in }S}} n-i,
$$
is a chain complex isomorphism.
\end{lemme}

\begin{proof}
Since $\psi$ is obviously an isomorphism of modules, we only need to check that the diagram is commutatif for any element of the basis.
Let $S$ and $S_0$ be enhanced states. In one hand, we have

\begin{eqnarray*}
\psi(d'(S))(S_0) & = & \sum_{\substack{v \textrm{ crossing in }D\\S' \textrm{ enhanced state}}} (-1)^{t^+_{S,v}+u(S')}[S',S]_v S'^*(S_0)\\
& = & \sum_{v \textrm{ crossing in }D}(-1)^{t^+_{S,v}+u(S_0)}[S_0,S]_v,
\end{eqnarray*}

and in the other

\begin{eqnarray*}
d^*(\psi(S))(S_0) & = & \sum_{\substack{v \textrm{ crossing in}D\\S' \textrm{ enhanced state}}} (-1)^{t^-_{S_0,v}+u(S)}[S_0,S']_v S^*(S')\\
& = & \sum_{v \textrm{ crossing in }D}(-1)^{t^-_{S_0,v}+u(S)}[S_0,S]_v.
\end{eqnarray*}

Moreover, if $[S_0,S]_v\neq 0$, then $S$ and $S_0$ differ only on $v$. We have thus $u(S_0)-u(S)=n-i$ for some $i$. On the other hand, we clearly have $t^+_{S,v}+t^-_{S_0,v}=t^+_{S,v}+t^-_{S,v}=n-i$. As a matter of fact, $(t^+_{S,v}+u(S_0))+(t^-_{S_0,v}+u(S))=2(n-i+u(S))$ and $t^+_{S,v}+u(S_0)$ and $t^-_{S_0,v}+u(S)$ have the same parity.
\end{proof}

\begin{lemme}
The map
$$
\phi : \vcenter{\hbox{\xymatrix{
\cdots \ar[r]^(.37){d} & \C_{i,j,k}(D) \ar[r]^(.45){d} \ar[d]^\phi & \C_{i-1,j,k}(D) \ar[r]^(.58){d} \ar[d]^\psi& \cdots\\
\cdots \ar[r]^(.32){d'} & \C'^{-i,-j,-k}(\overline{D}) \ar[r]^(.47){d'} & \C'^{-i+1,-j,-k}(\overline{D}) \ar[r]^(.66){d'} & \cdots,\\
}}}
$$
defined on any enhanced state $S$ by inversing the labels of all circles, is a chain complex isomorphism.
\end{lemme}

\begin{proof}
First, note that the underlying state of $S$ can be seen as a state of $\overline{D}$ with opposed smoothing at every crossing. Moreover, since $w(\overline{D})=-w(D)$ and since the map $\phi$ inverses all the labels, this map is well defined as a graded module isomorphism.Thus, we only need to check that the diagram is commutatif for any enhanced state $S$ of the usual basis.
Indeed, we have
$$
\phi(d(S)) = \sum_{\substack{v \textrm{ crossing in }D\\S' \textrm{ enhanced state of }D}} (-1)^{t^-_{S,v}}[S,S']_v \phi(S'),
$$
as well as
\begin{eqnarray*}
d'(\phi(S)) & = & \sum_{\substack{v \textrm{ crossing in }\overline{D}\\S' \textrm{ enhanced state of }\overline{D}}} (-1)^{t^+_{\phi(S),v}}[S',\phi(S)]_v S'\\
& = & \sum_{\substack{v \textrm{ crossing in }\overline{D}\\\phi(S') \textrm{ enhanced state of }\overline{D}}} (-1)^{t^+_{\phi(S),v}}[\phi(S'),\phi(S)]_v \phi(S'),
\end{eqnarray*}
since $\phi$ is a bijection between enhanced states of $D$ and $\overline{D}$. But, by definition, $t^-_{S,v}=t^+_{\phi(S),v}$ and by construction $[\phi(S'),\phi(S)]_v = [S,S']_v$.
\end{proof}

Now, the proposition is obtained by passing to homology the precedent lemnata.

\subsection{Invariance under star--like Reidemeister moves of type II}

Here, we will follow the proofs given in \cite{Audoux2}, details are done there. Modules are depicted by sets of generators and sets of generators by diagrams partially smoothed. A given set is obtained by considering all the choices for smoothing the unsmoothed crossings and labeling the unlabeled circles. This description is also used for formulas, the choice is then consistent all over the equality. Partial differentials associated to a crossing $c$, restricted to a given subcomplex $C$, are denoted with a sequence of indices which correspond, increasingly, to the crossings which are smoothed in the same way for all generators of $C$. A star is assigned to $c$, otherwise, the index corresponds to the nature of the common smoothing.\\

We introduce the operators $.\{n\}$ (resp. $.[n]$) which is the global uplifting of the grading (resp. homological) degree by $n$. The key ingredients for upcoming proofs are the two following propositions~:

\begin{defi}
Let $\func{f}{\C}{\C'}$ be a morphism of chain complexes. Then the cone of $f$ is the chain complex $\C''$ defined by $\C'' = \C \oplus \C'[-1]$ and $d_{\C''}=(d_\C + f)\oplus (-d_{\C'})$.
\end{defi}

\begin{prop}\label{util1}
Let $D$ be a diagram of $K$ and $v$ a crossing of $D$. Let $D_0$ and $D_1$ be the diagrams obtained from $D$ by smoothing $v$ respectively in the $A$--fashion and in the $A^{-1}$ one. Then $d_v$ defines a grading preserving map of chain complexes from $\C(D_0)$ to $\C(D_1)\{-1\}$. Moreover, $\C(D)$ is the cone of this map.
\end{prop}

\begin{prop}\label{util2}
The cone of an isomorphism between two chain complexes is acyclic.
\end{prop}

\vspace{.7cm}

Now consider two diagrams which differ only locally by a star--like Reidemeister move of type II.
$$
\begin{array}{ccc}
D & \hspace{2cm} & D'\\
\dessin{1.1cm}{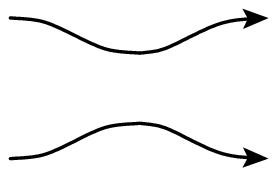} & & \dessin{1.1cm}{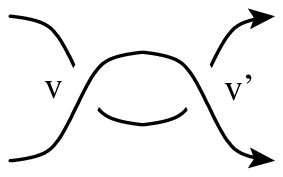}.
\end{array}
$$
We order the crossings of $D$ and $D'$ in the same fashion, letting $v$ and $v'$ be the last two ones of $D'$. Then we consider the following four diagrams which are partial smoothings of $D'$~:
$$
\begin{array}{cccl}
D'_{00} &  \hspace{1.8cm} & D'_{10}&\\
\dessin{1cm}{II1} &  & \dessin{1cm}{II2}&\\[5mm]
D'_{01} &  & D'_{11}&\\
\dessin{1cm}{II3} & & \dessin{1cm}{II4}&,\\
\end{array}
$$
and two maps~:
$$
\begin{array}{c}
  \func{d_{0\star}}{\left\{\dessin{.8cm}{II1}\right\}}{\left\{\dessin{.8cm}{II3}\right\}},\\[.5cm]
  \func{d_{\star  1}}{\left\{\dessin{.8cm}{II3}\right\}}{\left\{\dessin{.8cm}{II4}\right\}}.
\end{array}
$$

\begin{lemme}
The morphisms $d_{0\star}$ and $d_{\star 1}$ are respectively injective and surjective.
\end{lemme}

\begin{lemme}\label{firstproof}
The cone of
$$
\xymatrix @!0 @C=3cm @R=1.7cm {
\cdots  &  Im(d^i_{\star 0}+ d^i_{0 \star}) \ar[l]_(.54){d^{i-1}_{D'}} & \cdots \ar[l]_(.35){d^i_{D'}}\\
\cdots  & \left\{\dessin{.8cm}{II1}\right\}  \ar[u]^(.55){d^i_{\star 0}+ d^i_{0 \star}} \ar[l]_(.6){d^{i-1}_{D'_{00}}}& \cdots \ar[l]_(.37){d^i_{D'_{00}}}\\
}
$$
is an acyclic graded subcomplex of $\C(D')$ denoted $\C_1$.
\end{lemme}

\begin{lemme}
The cone of
$$
\xymatrix @!0 @C=3cm @R=1.7cm {
\cdots &  \left\{\dessin{.8cm}{II4}\right\} \ar[l]_(.58){d^{i-1}_{D'}}  & \cdots \ar[l]_(.38){d^i_{D'}}\\
\cdots   & \left\{\dessin{.8cm}{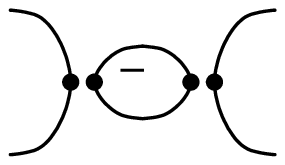}\right\}  \ar[u]^{d^i_{\star 1}} \ar[l]_(.58){d^{i-1}_{D'_{01}}} & \cdots \ar[l]_(.38){d^i_{D'_{01}}}
}
$$
is an acyclic graded subcomplex of $\C(D')$ denoted $\C_2$.
\end{lemme}

\begin{defi}
We define the map
$$
\begin{array}{rccc}
.\otimes v_- : & \left\{\dessin{.8cm}{II4}\right\} & \longrightarrow & \left\{\dessin{.85cm}{II4}\right\}\\[4mm]
& \dessin{.65cm}{II4} & \mapsto & -\dessin{.65cm}{II3m}.
\end{array}
$$
\end{defi}

\begin{remarque} $\ $
The map $.\otimes v_-$ is grading--preserving and is a right inverse for $d_{\star 1}$.
\end{remarque}

\begin{lemme}
The complex
$$
\xymatrix @!0 @C=6cm @R=2cm {
\cdots  &  \left\{\dessin{.8cm}{II2}-\big( d_{1 \star}(\dessin{.8cm}{II2}) \big)\otimes v_- \right\} \ar[l]_(.67){d^{i-1}_{D'}} & \cdots \ar[l]_(.28){d^i_{D'}}}
$$
is a graded subcomplex of $\C(D')$ denoted $\C_3$, isomorphic to $\C(D)$ via $\psi_{II}$ defined by
$$
\psi_{II}\left(\dessin{.7cm}{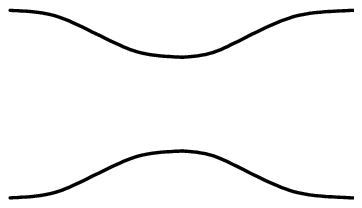} \right)=\dessin{.8cm}{II2}-\big( d_{1 \star}(\dessin{.8cm}{II2}) \big)\otimes v_-.
$$
\end{lemme}

\begin{lemme}\label{lastlemme}
$\C (D') \simeq \C_1 \oplus \C_2 \oplus \C_3$.
\end{lemme}

Finally, we conclude by passing to homolgy in Lemma \ref{lastlemme}.

\subsection{Invariance under star--like Reidemeister moves of type III}

Since the mirror image of a move of type $III_g$ is a move of type $III_h$ and as the Corollary \ref{mirroir} holds, it is sufficient to prove the invariance under Reidemeister moves of type $III_h$.\\

Consider two diagrams which differ only locally by the following star--like Reidemeister move of type $III_h$.
$$
\begin{array}{ccc}
D & \hspace{2cm} & D'\\
\dessin{1.375cm}{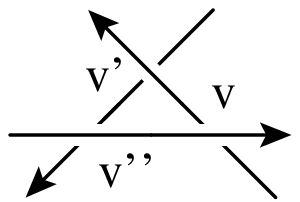} & & \dessin{1.375cm}{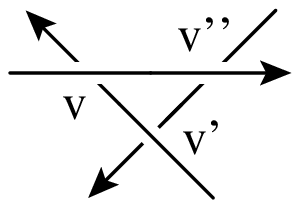}.
\end{array}
$$
We order the crossings of $D$ and $D'$ in the same fashion, letting the three crossings involved in the Reidemeister move be the last three ones. We can now consider the following ten partial smoothings of $D$ or $D'$~:
$$
\begin{array}{ccc}
\begin{array}{ccc}
D_{000}& \hspace{.8cm} & D_{001}\\
\dessin{1.25cm}{III2s2} &  & \dessin{1.25cm}{III2s1}
\end{array}
&  \hspace{1cm} &
\begin{array}{ccc}
D'_{000}&  \hspace{.8cm} & D'_{001}\\
\dessin{1.25cm}{III12} & & \dessin{1.25cm}{III11}
\end{array}\\[10mm]
\begin{array}{ccc}
D_{010}&  \hspace{.8cm} & D_{011}\\
\dessin{1.25cm}{III3s3} & & \dessin{1.25cm}{III2s4}
\end{array}
&  \hspace{1cm} &
\begin{array}{ccc}
D'_{010}&  \hspace{.8cm} & D'_{011}\\
\dessin{1.25cm}{III33} & & \dessin{1.25cm}{III13}
\end{array}\\[10mm]
D_{1\bullet \bullet}&  \hspace{1.2cm} & D'_{1\bullet \bullet}\\
\dessin{1.25cm}{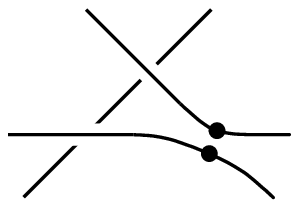} & & \dessin{1.25cm}{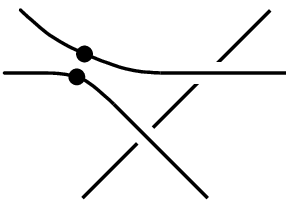}.
\end{array}
$$
One can easily check that $D_{1\bullet \bullet}\simeq D'_{1\bullet \bullet}$ as chains complexes.\\
With respect to $D$, we also consider two maps~:
$$
\begin{array}{cccccl}
d_{00\star} &:& \left\{\dessin{1cm}{III2s2}\right\} &\longrightarrow& \left\{\dessin{1cm}{III2s1}\right\}&\\[.5cm]
d_{0\star 1} &:&  \left\{\dessin{1cm}{III2s1}\right\}&\longrightarrow& \left\{\dessin{1cm}{III2s4}\right\}&.
\end{array}
$$

\begin{lemme}
The maps $d_{00\star}$ and $d_{0\star 1}$ are, respectively injective and surjective.
\end{lemme}

Then, the rest follows as in the case of the Reidemeister move of type II.

\begin{lemme}
The cone of
$$
\xymatrix @!0 @C=4cm @R=2cm {
\cdots &  Im(d^i_{\star 00}+ d^i_{0 \star 0} + d^i_{00\star}) \ar[l]_(.7){d^{i-1}_D} & \cdots \ar[l]_(.22){d^i_D}\\
\cdots & \Z\left\{\dessin{1cm}{III2s2}\right\} \ar[u]^{d^i_{\star 00}+ d^i_{0 \star 0}+d^i_{00\star}} \ar[l]_(.55){d^{i-1}_{D_{000}}} & \cdots \ar[l]_(.41){d^i_{D_{000}}}\\
}
$$
is an acyclic graded subcomplex of $\C(D)$ denoted $\C_1$.
\end{lemme}

\begin{lemme}
The cone of
$$
\xymatrix @!0 @C=4,2cm @R=2cm {
\cdots & \left\{\dessin{1cm}{III2s3}+\dessin{1cm}{III2s4}\right\} \{-1\} \ar[l]_(.67){d^{i-1}_D} & \cdots \ar[l]_(.26){d^i_D} \\
\cdots  & \left\{\dessin{1cm}{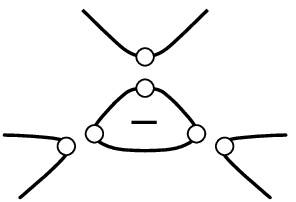}\right\}  \ar[u]^{d^i_{\star 01}+d^i_{0\star 1}} \ar[l]_(.54){d^{i-1}_{D_{001}}} & \cdots \ar[l]_(.41){d^i_{D_{001}}}
}
$$
is an acyclic graded subcomplex of $\C(D')$ denoted $\C_2$.
\end{lemme}

\begin{defi}
We define the map
$$
\begin{array}{rccc}
.\otimes v_- : & \left\{\dessin{1cm}{III2s4}\right\} & \longrightarrow & \left\{\dessin{1cm}{III2s1}\right\}\\[4mm]
& \dessin{.8125cm}{III2s4} & \mapsto & -\dessin{.8125cm}{III2s1m}.
\end{array}
$$
\end{defi}

\begin{remarque} $\ $
The map $.\otimes v_-$ is grading--preserving and is a right inverse for $d_{0\star 1}$.
\end{remarque}

\begin{lemme}
The complex
$$
\xymatrix @!0 @C=6.5cm @R=2cm {
\cdots &  \left\{\dessin{1cm}{III3s3}-\big( d_{01 \star}(\dessin{1cm}{III3s3}) \big)\otimes v_- \right\}\oplus \left\{\dessin{1cm}{IIIspp}\right\} \ar[l]_(.73){d^{i-1}_D} & \cdots \ar[l]_(.23){d^i_D} }
$$
is a graded subcomplex of $\C(D)$ denoted $\C_3$.
\end{lemme}

\begin{lemme}
$\C(D) \simeq \C_1 \oplus \C_2 \oplus \C_3$.
\end{lemme}

With the same reasoning on $D'$, we get~:
\begin{lemme}
$\C(D') \simeq \C'_1 \oplus \C'_2 \oplus \C'_3$ with $\C'_1$ and $\C'_2$ acyclic and
$$
\C'_3 = \left\{\dessin{1cm}{III33}-\big( d_{01 \star}(\dessin{1cm}{III33}) \big)\otimes v_- \right\}\oplus \left\{\dessin{1cm}{IIIpp}\right\}.
$$
\end{lemme}

Finally, we conclude thanks to the following lemma~:

\begin{lemme}
The chain complexes $\C_3$ and $\C'_3$ are isomorphic via $\psi_{III}$ defined on $\left\{\dessin{1cm}{III3s3}\right\}$ by
$$
\psi\left(\dessin{1cm}{III3s3}-\big( d_{01 \star}(\dessin{1cm}{III3s3}) \big)\otimes v_- \right)=\dessin{1cm}{III33}-\big( d_{01 \star}(\dessin{1cm}{III33}) \big)\otimes v_-
$$
and on $\left\{\dessin{1cm}{IIIspp}\right\}$ in the obvious way.
\end{lemme}

\subsection{Invariance under star--like Reidemeister moves of type I}

Consider the three following diagrams which differ only locally by a Reidemeister move of type I and the adding of a circle.
$$
\begin{array}{ccccc}
D' & \hspace{1.5cm} & D & \hspace{1.5cm} & D''\\
\dessin{1.21cm}{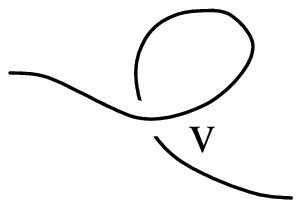} & & \dessin{1.2cm}{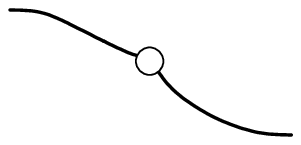} & & \dessin{1.21cm}{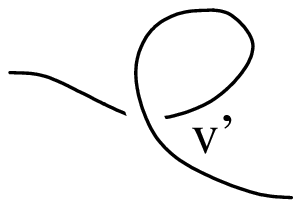}.
\end{array}
$$
We order the crossings of $D$, $D'$ and $D''$ in the same fashion, letting $v$ and $v'$ be the last ones of, respectively, $D'$ and $D''$. Then we consider the following diagrams which are partial smoothings of $D'$ and $D''$~:
$$
\begin{array}{cccl}
D'_{0} & \hspace{1.8cm} & D''_{0}&\\
\dessin{1.1cm}{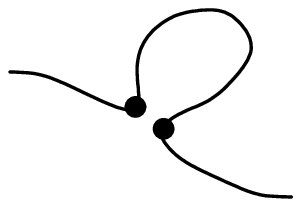} & & \dessin{1.1cm}{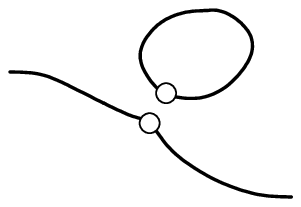}&\\[5mm]
D'_{1} & & D''_{1}&\\
\dessin{1.1cm}{RI11} & & \dessin{1.1cm}{RI10}&,\\
\end{array}
$$
and two maps~:
$$
\begin{array}{cccccl}
d'_{\star} &:& \left\{\dessin{.88cm}{RI10}\right\} &\longrightarrow& \left\{\dessin{.88cm}{RI11}\right\}&\\[.5cm]
d''_{\star} &:&  \left\{\dessin{.88cm}{RI11}\right\}&\longrightarrow& \left\{\dessin{.88cm}{RI10}\right\}&.
\end{array}
$$

\begin{lemme}
The morphisms $d'_{\star}$ and $d''_{\star}$ are respectively injective and surjective.
\end{lemme}

\begin{lemme}\label{secondproof}
The cone of
$$
\xymatrix @!0 @C=3cm @R=2.5cm {
\cdots  &  \left\{\dessin{.88cm}{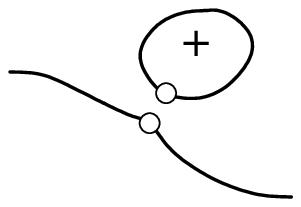}\right\} \ar[l]_(.68){d^{i-1}_{D'}} & \cdots \ar[l]_(.26){d^i_{D'}}\\
\cdots  & \left\{\dessin{.88cm}{RI10}\right\}  \ar[u]^{d'_{\star}} \ar[l]_(.6){d^{i-1}_{D'_0}}& \cdots \ar[l]_(.33){d^i_{D'_0}}\\
}
$$
is an acyclic graded subcomplex of $\C(D')$ denoted $\C'_1$.
\end{lemme}

\begin{lemme}
The cone of
$$
\xymatrix @!0 @C=3cm @R=2.5cm {
\cdots & \left\{\dessin{.88cm}{RI10}\right\} \ar[l]_(.65){d^{i-1}_{D''}}  & \cdots \ar[l]_(.27){d^i_{D'}}\\
\cdots   & \left\{\dessin{.88cm}{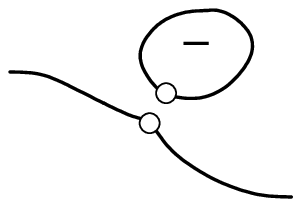}\right\}  \ar[u]^{d''_{\star}} \ar[l]_(.61){d^{i-1}_{D''_0}} & \cdots \ar[l]_(.35){d^i_{D''_0}}
}
$$
is an acyclic graded subcomplex of $\C(D'')$ denoted $\C''_1$.
\end{lemme}

\begin{lemme}
The complexes
$$
\xymatrix @!0 @C=4cm @R=2cm {
\cdots  &  \left\{\dessin{.88cm}{RI11m}\right\} \ar[l]_(.55){d^{i-1}_{D'}} & \cdots \ar[l]_(.4){d^i_{D'}}}
$$
and
$$
\xymatrix @!0 @C=4,5cm @R=2cm {
\cdots  &  \left\{\dessin{.88cm}{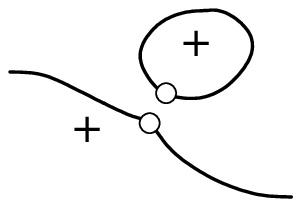},\dessin{.88cm}{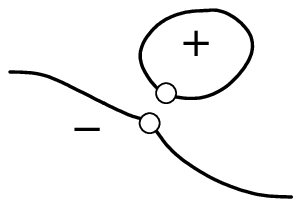}-\dessin{.88cm}{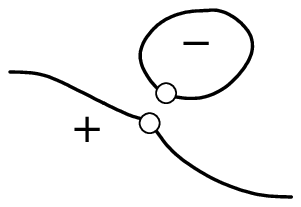}\right\} \ar[l]_(.75){d^{i-1}_{D''}} & \cdots \ar[l]_(.2){d^i_{D''}}}
$$
are graded subcomplexes of respectively $\C(D')$ and $\C(D'')$ denoted $\C'_2$ and $\C''_2$. Both of them are isomorphic to $\C(D)$.
\end{lemme}

\begin{lemme}\label{lastlemme2}
$\C (D') \simeq \C'_1 \oplus \C'_2$ and $ \C(D'')\simeq\C''_1 \oplus \C''_2$.
\end{lemme}

\begin{lemme}
The adding of two adjacent Seifert points on a diagram do not change its chain complex.
\end{lemme}

After having passed to homology the precedent lemma, one can easily achieve the invariance under star--like Reidemeister moves of type I.\\

\section{A picture for the invariant}
\label{sec:Picture}

In a somehow similar way than Bar-Natan (\cite{BarNatan}), we can give a picture for the above construction.
Details are done in a forthcoming paper (\cite{Audoux3}) but Figure 10 gives the main ingredients.\\
Roughly speaking, it corresponds to a cube of resolution of which vertices are obtained by associating a dot to every $h$-circle and a circle for each $d$-one.
Edges are then decorated in a standard way by surfaces with pulleys.\\

\begin{figure}[p]

  \subfigure[Generators]{
    \framebox[14cm]{
      \begin{minipage}{1.0\linewidth}
        \begin{center}
          $\vcenter{\hbox{\includegraphics[height=1.5cm]{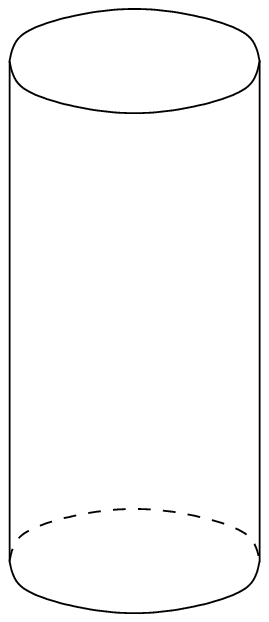}}}
          \hspace{.5cm}
          \vcenter{\hbox{\includegraphics[height=1.5cm]{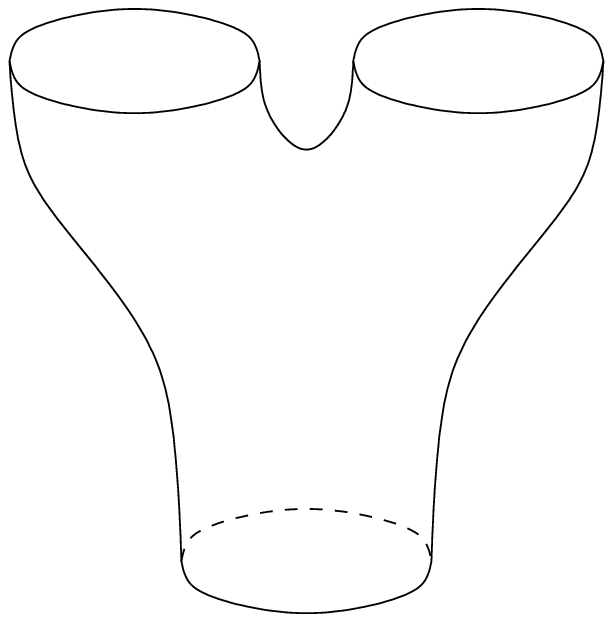}}}
          \hspace{.5cm}
          \vcenter{\hbox{\includegraphics[height=1.5cm]{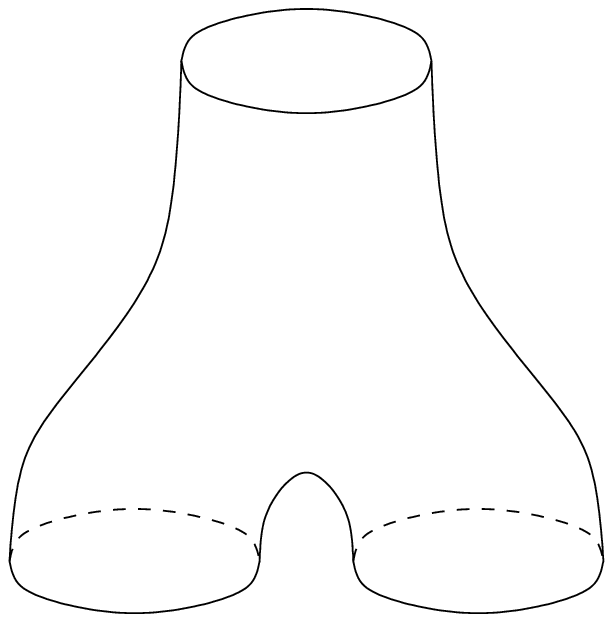}}}
          \hspace{.5cm}
          \vcenter{\hbox{\includegraphics[height=1.5cm]{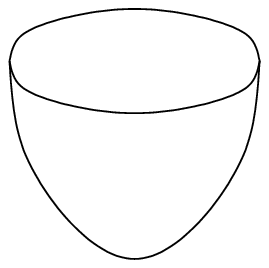}}}
          \hspace{.5cm}
          \vcenter{\hbox{\includegraphics[height=1.5cm]{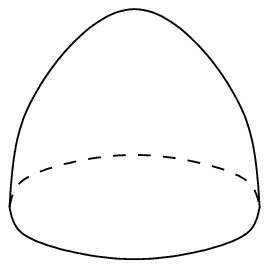}}}
          \hspace{.5cm}
          \vcenter{\hbox{\includegraphics[height=1.5cm]{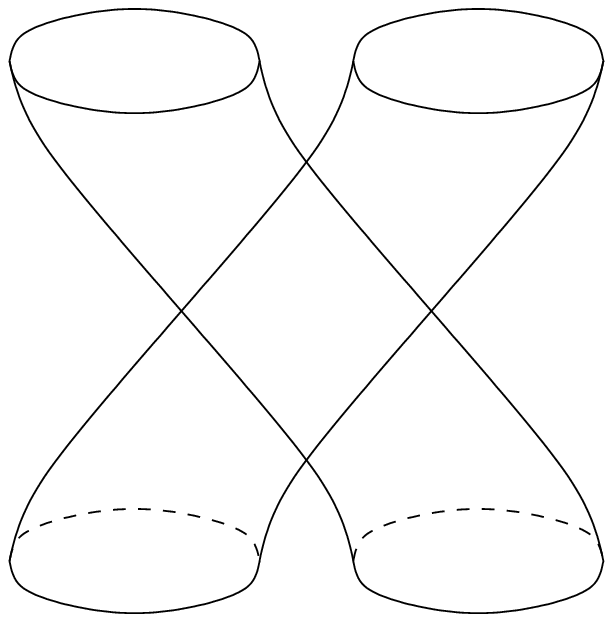}}}$\\[.2cm]
          $2$--dimensional generators
        \end{center}
        
        \vspace{.6cm}
        
        \begin{center}
          $\vcenter{\hbox{\includegraphics[height=1.5cm]{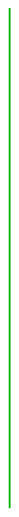}}}
          \hspace{.5cm}
          \vcenter{\hbox{\includegraphics[height=1.5cm]{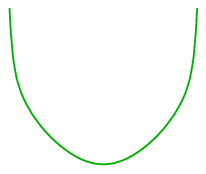}}}
          \hspace{.5cm}
          \vcenter{\hbox{\includegraphics[height=1.5cm]{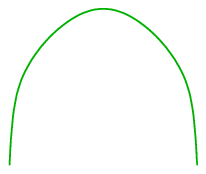}}}
          \hspace{.5cm}
          \vcenter{\hbox{\includegraphics[height=1.5cm]{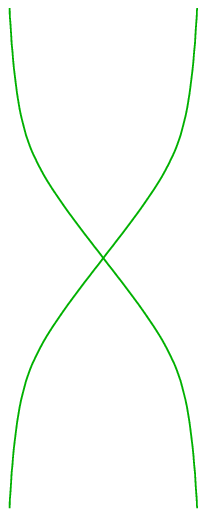}}}$\\[.2cm]
          $1$--dimensional generators
        \end{center}
        
        \vspace{.6cm}
        
        \begin{center}
          $\vcenter{\hbox{\includegraphics[height=1.5cm]{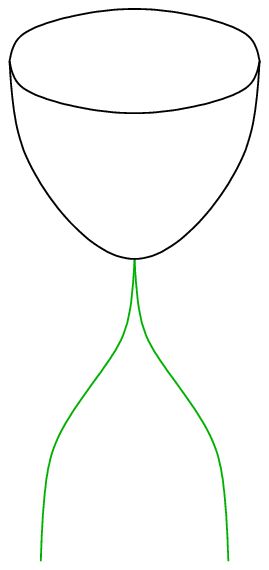}}}
          \hspace{.5cm}
          \vcenter{\hbox{\includegraphics[height=1.5cm]{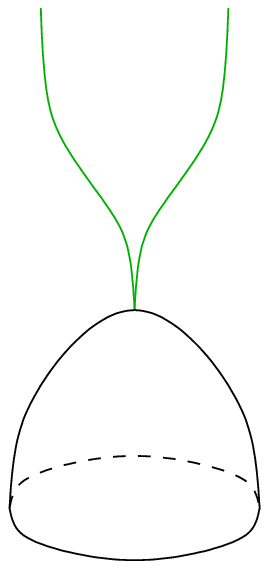}}}
          \hspace{.5cm}
          \vcenter{\hbox{\includegraphics[height=1.5cm]{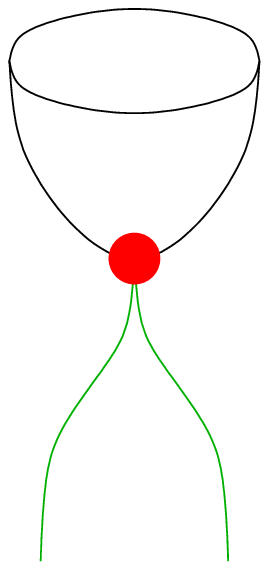}}}
          \hspace{.5cm}
          \vcenter{\hbox{\includegraphics[height=1.5cm]{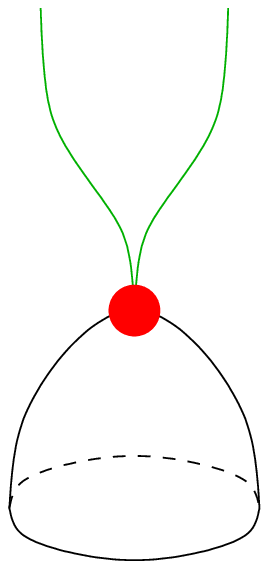}}}
          \hspace{.5cm}
          \vcenter{\hbox{\includegraphics[height=1.5cm]{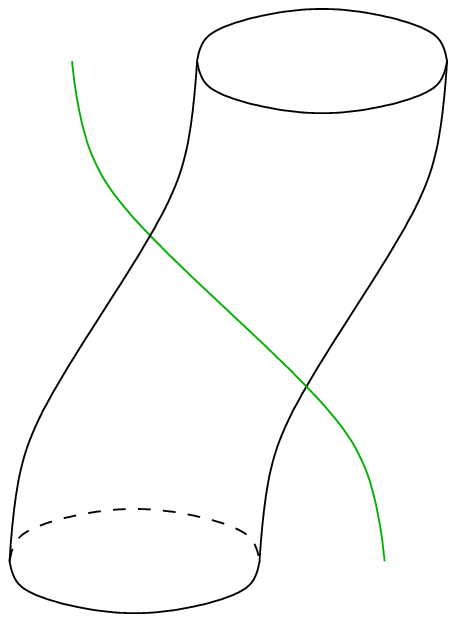}}}
          \hspace{.5cm}
          \vcenter{\hbox{\includegraphics[height=1.5cm]{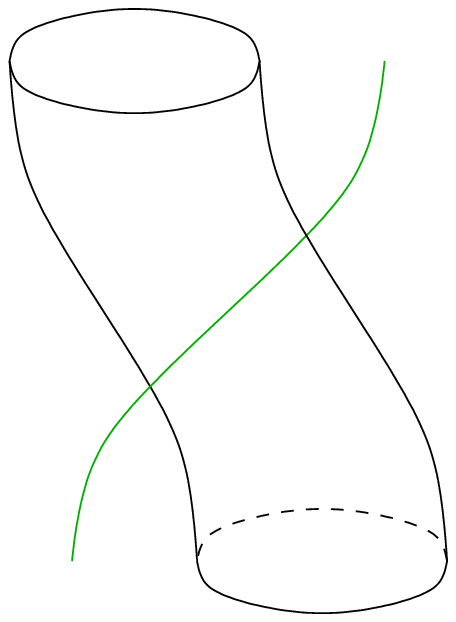}}}$\\[.2cm]
          mixed generators
        \end{center}
      \end{minipage}
    }
    \label{fig:generators}} 
  
  \subfigure[Relations]{
    \framebox[14cm]{
      \begin{minipage}{1.0\linewidth}
        \begin{center}
          $B:\hspace{.3cm} \dessin{1.5cm}{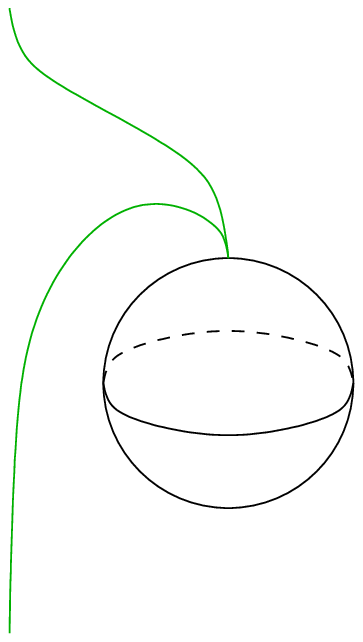}\ =\ \dessin{1.5cm}{1gen1}
          \hspace{2cm}
          R:\hspace{.3cm}  \dessin{1.5cm}{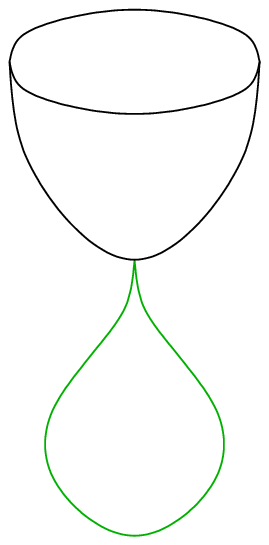}\ =\ \dessin{1.5cm}{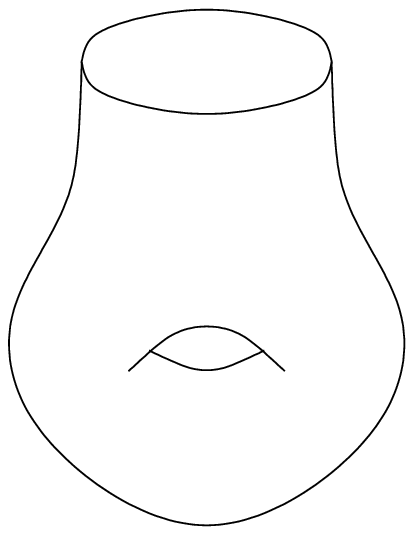}$\\[.2cm]
          Adding/deleting bubles and ropes
        \end{center}
        
        \vspace{.6cm}
        
        \begin{center}
          $P:\hspace{.3cm} \dessin{1.5cm}{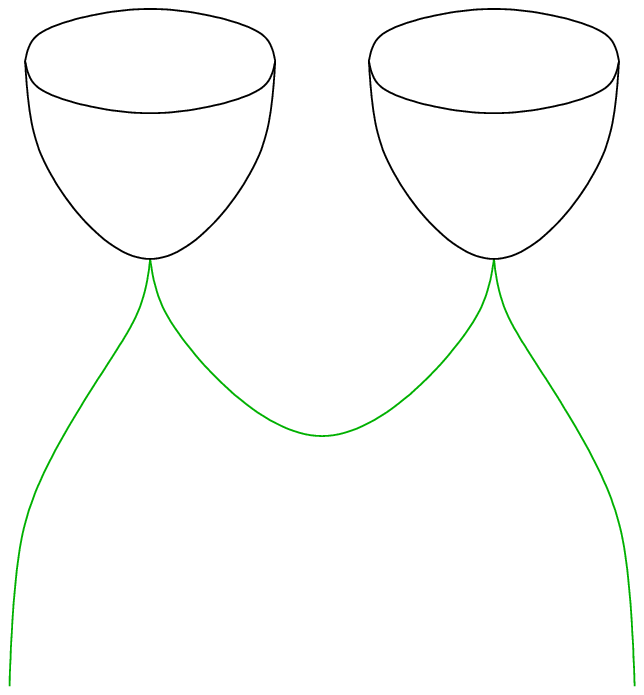}\ =\ \dessin{1.5cm}{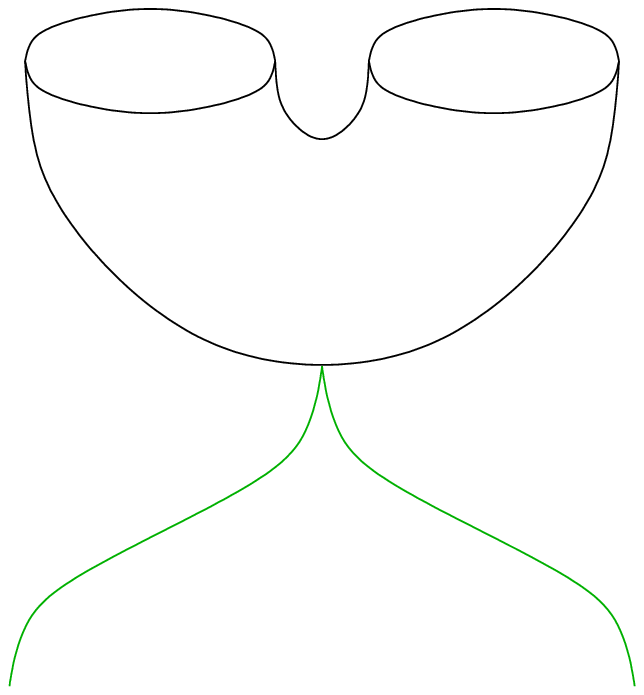}
          \hspace{2cm}
          P':\hspace{.3cm} \dessin{1.5cm}{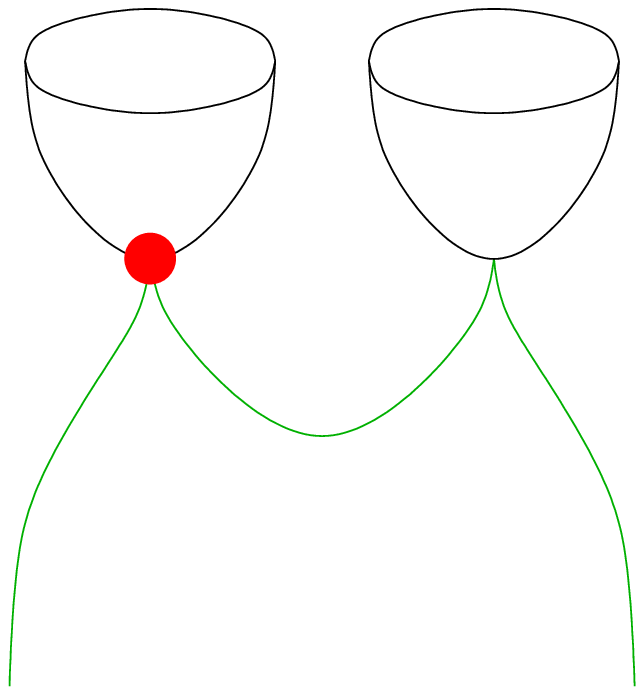}\ =\ \dessin{1.5cm}{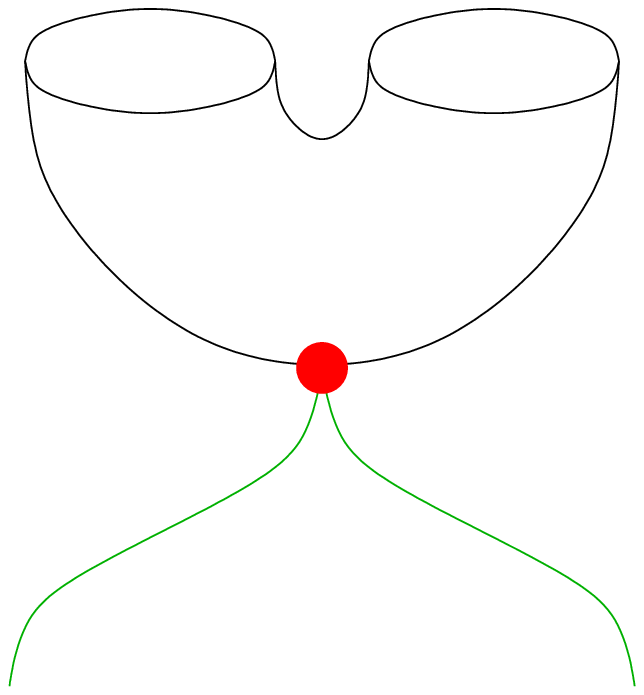}$\\[.2cm]
          Pulling ropes
        \end{center}
        
        \vspace{.6cm}
        
        \begin{center}
          $L:\hspace{.3cm} \dessin{1.5cm}{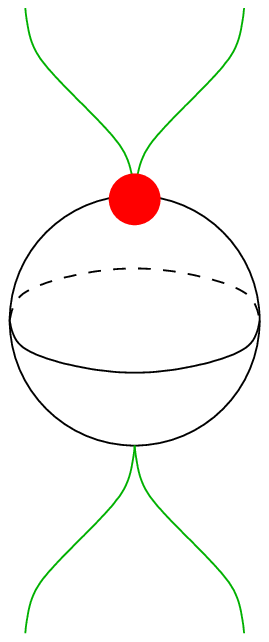}\ =\ \dessin{1.5cm}{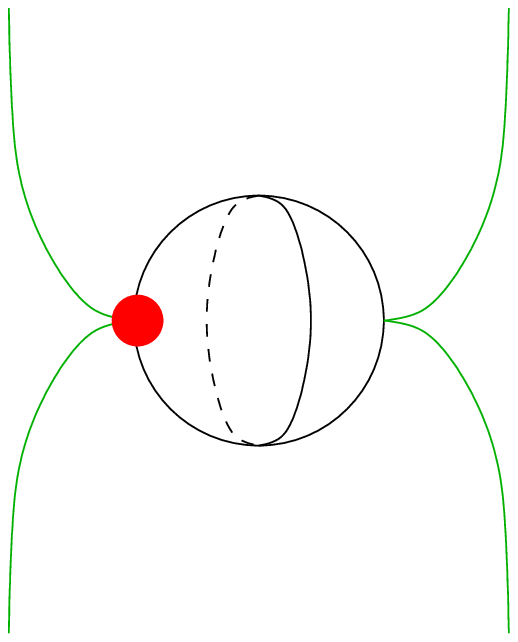}
          \hspace{2cm}
          D: \hspace{.3cm} \dessin{1.5cm}{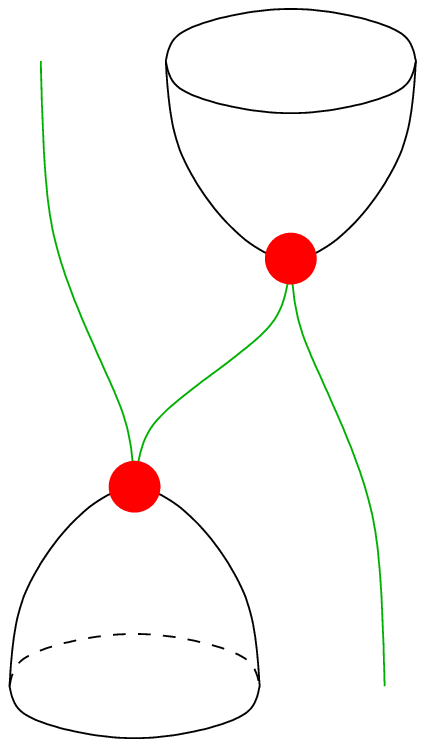}\ =\ \dessin{1.5cm}{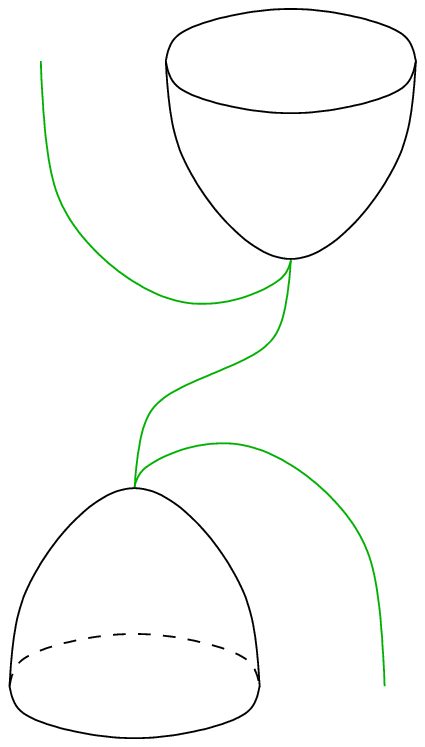}$\\[.2cm]
          $2P: \hspace{.3cm} \dessin{1.5cm}{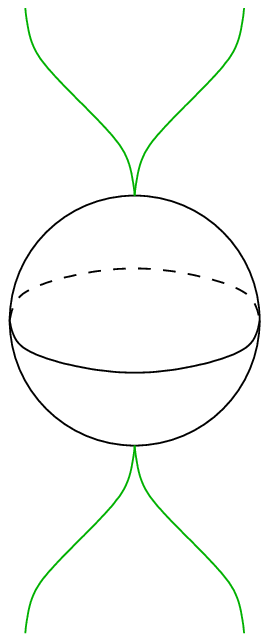}\ =\ 0 \ =\  \dessin{1.5cm}{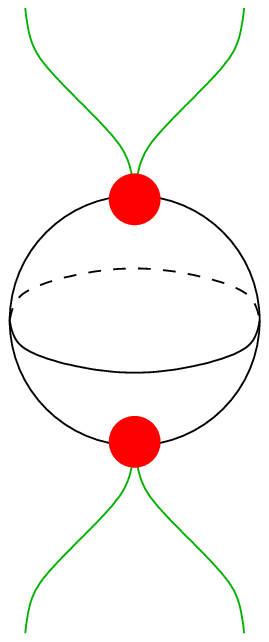}  \hspace{.3cm} :2P'$\\[.2cm]
          Involving two pulleys
        \end{center}
      \end{minipage}
    }
    
    \label{fig:relations}}
  
  \label{fig:DefinitionSurfaces}
  \caption{Definition of surfaces with pulleys}
\end{figure}

Surfaces with pulley are abstract surfaces composed by $1$ and $2$--dimensional pieces that can be connected together on a finite number of points of their interior. These points are called pulleys and we distinguish two kind of them (pictured with or without a dot).
In Figure \ref{fig:generators}, we give a set of generators.
Any two surfaces can be composed horizontally by disjoint union.
Moreover, if the top boundary of a surface corresponds to the bottom one of another surface, then they can be composed vertically by gluing them along this common boundary.

\begin{figure}[p]
  \centering
$$
\begin{array}{c}L :\\[.7cm] \ \end{array} \dessin{2.2cm}{K} \hspace{9cm}
$$

\vspace{-1cm}

\rotatebox{340}{$\leadsto$} \hspace{4cm}

\vspace{-2.5cm}

$$
\hspace{4cm}
\xymatrix{
& \dessin{1.3cm}{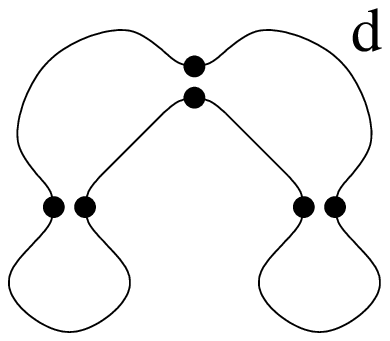} \ar[r] \ar[dr]&  \dessin{1.3cm}{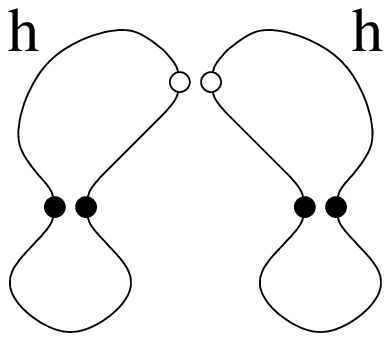} \ar[dr]& \\
\dessin{1.3cm}{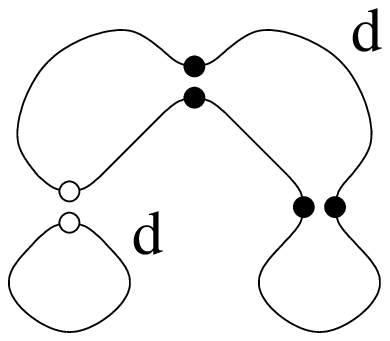} \ar[ur] \ar[r] \ar[dr] &  \dessin{1.3cm}{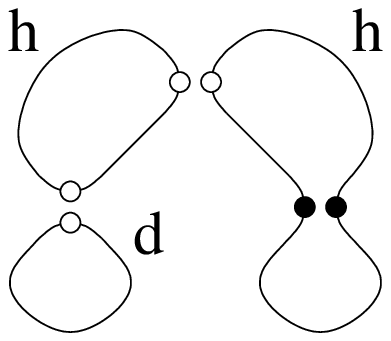} \ar[ur]|!{[u];[r]}\hole \ar[dr]|!{[d];[r]}\hole &  \dessin{1.3cm}{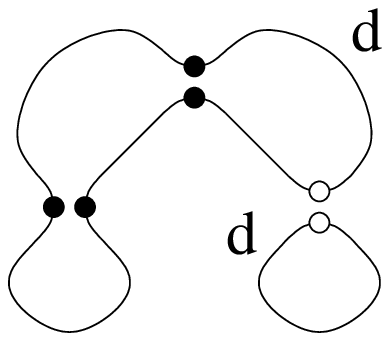} \ar[r]& \dessin{1.3cm}{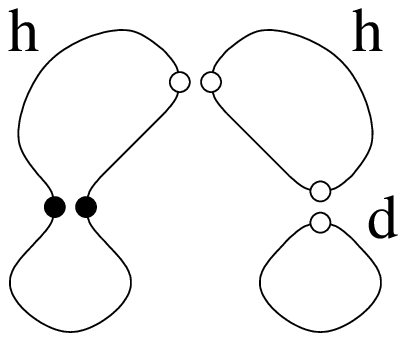}\\
& \dessin{1.3cm}{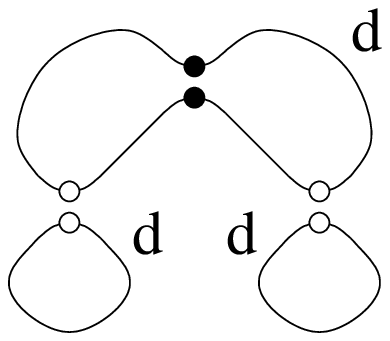} \ar[ur] \ar[r]&  \dessin{1.3cm}{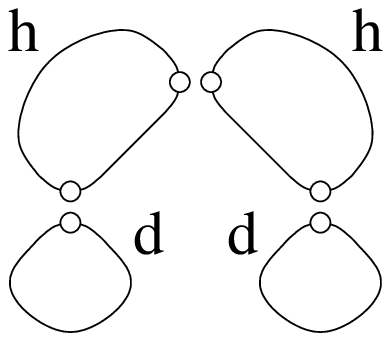} \ar[ur]&}
$$

\vspace{-.3cm}

\hspace{1.5cm} \rotatebox{225}{$\leadsto$}

\vspace{-.2cm}

$$
\xymatrix@!0@C=2cm@R=4cm{
&& \dessinH{1.1cm}{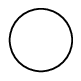} \ar[rr]^{\dessinH{1.1cm}{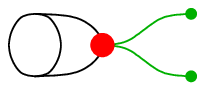}} \ar[drr]^(.3){\dessinH{1.1cm}{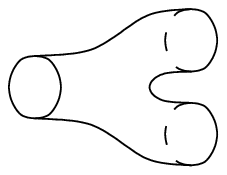}} \ar@{.}[d]|{\oplus}&&  \dessinH{1.1cm}{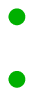} \ar[drr]^(.38){\dessinH{1.1cm}{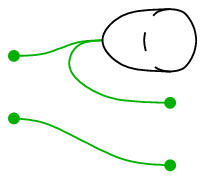}}\ar@{.}[d]|{\oplus}&& \\
\dessinH{1.1cm}{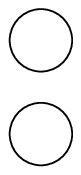} \ar[urr]^(.62){\dessinH{1.1cm}{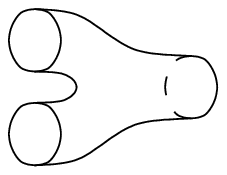}} \ar[rr]^(.48){\dessinH{1.1cm}{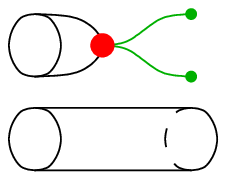}} \ar[drr]_(.62){\dessinH{1.1cm}{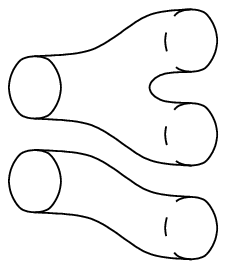}} \ar@{.>}[dd] & \ar@{}[dd]|(.5){}="H1" &  \dessinH{1.1cm}{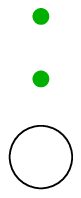} \ar[urr]^(.35){\dessinH{1.1cm}{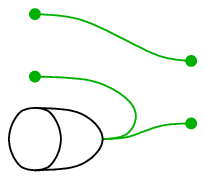}}|!{[u];[rr]}\hole \ar[drr]^(.3){\dessinH{1.1cm}{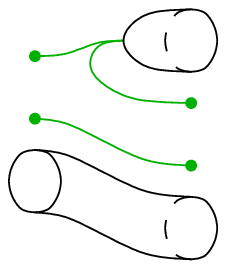}}|!{[d];[rr]}\hole \ar@{.}[d]|{\oplus}&&  \dessinH{1.1cm}{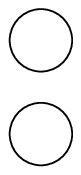} \ar[rr]^(.48){\dessinH{1.1cm}{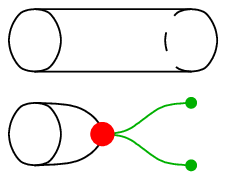}}\ar@{.}[d]|{\oplus}& \ar@{}[dd]|(.5){}="H3"& \dessinH{1.1cm}{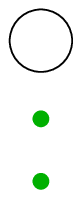}\ar@{.>}[dd]\\
&& \dessinH{1.1cm}{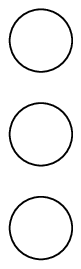} \ar[urr]^(.35){\dessinH{1.1cm}{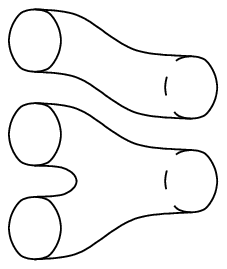}} \ar[rr]_{\dessinH{1.1cm}{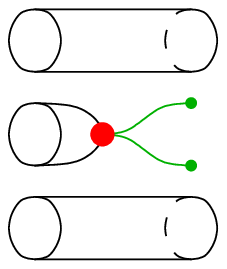}} \ar@{.>}[d]& \ar@{}[d]|(.42){}="H2"& \dessinH{1.1cm}{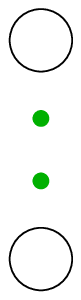} \ar[urr]_(.38){\dessinH{1.1cm}{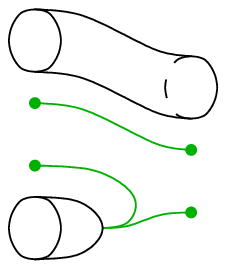}} \ar@{.>}[d]
&&\\
Ch_2 \ar[rr]^\partial="B1"  && Ch_1 \ar[rr]^\partial="B2" && Ch_0 \ar[rr]^\partial="B3" && Ch_{-1} \ar@{-->}"H1";"B1"  \ar@{-->}"H2";"B2"  \ar@{-->}"H3";"B3"}
$$

  \label{fig:Cube}
  \caption{Construction of the invariant: {\footnotesize for the sake of clarity, surfaces morphisms are depicted from the left to the right.}}
\end{figure}

Besides Bar-Natan $S$ and $4T$ relations, and besides the usual and obvious ones, some extra relations are also required.
They are given in Figure \ref{fig:relations}.
For the sake of clarity, ropes and surfaces are depicted as embedded in $[0,1]\times\R^2$, but the reader should keep in mind that every object can cross each other.\\

The construction given in this paper can be recovered by applying a suitable fonctor.
By changing the definition of $h$-circles, which are send to dots, the above picture can straighforwardly be adapted to the braid-like case and to the $I$-bundle case.

\nocite{*}
\bibliographystyle{amsalpha}
\bibliography{StarLikeJones}

\providecommand{\bysame}{\leavevmode\hbox to3em{\hrulefill}\thinspace}
\providecommand{\MR}{\relax\ifhmode\unskip\space\fi MR }
\providecommand{\MRhref}[2]{%
  \href{http://www.ams.org/mathscinet-getitem?mr=#1}{#2}
}
\providecommand{\href}[2]{#2}
\begin{thebibliography}{HGR05}

\bibitem[AF05]{Audoux2}
Benjamin Audoux and Thomas Fiedler, \emph{A {J}ones polynomial for braid-like
  isotopies of oriented links and its categorification}, Algebr. Geom. Topol.
  \textbf{5} (2005), 1535--1553 (electronic). \MR{MR2186108 (2006h:57008)}

\bibitem[APS04]{Asaeda}
Marta~M. Asaeda, J{\'o}zef~H. Przytycki, and Adam~S. Sikora,
  \emph{Categorification of the {K}auffman bracket skein module of
  {$I$}-bundles over surfaces}, Algebr. Geom. Topol. \textbf{4} (2004),
  1177--1210 (electronic). \MR{MR2113902 (2006a:57010)}

\bibitem[Aud07]{Audoux3}
Benjamin Audoux, \emph{Category of surfaces with pulleys and khovanov
  homology}, in preparation, 2007.

\bibitem[BN05]{BarNatan}
Dror Bar-Natan, \emph{Khovanov's homology for tangles and cobordisms}, Geom.
  Topol. \textbf{9} (2005), 1443--1499 (electronic). \MR{MR2174270
  (2006g:57017)}

\bibitem[BZ03]{Burde}
Gerhard Burde and Heiner Zieschang, \emph{Knots}, second ed., de Gruyter
  Studies in Mathematics, vol.~5, Walter de Gruyter \& Co., Berlin, 2003.
  \MR{MR1959408 (2003m:57005)}

\bibitem[Fie01]{Fiedler2}
Thomas Fiedler, \emph{Gauss diagram invariants for knots and links},
  Mathematics and its Applications, vol. 532, Kluwer Academic Publishers,
  Dordrecht, 2001. \MR{MR1948012 (2003m:57031)}

\bibitem[HGR05]{Helme}
Laure Helme-Guizon and Yongwu Rong, \emph{A categorification for the chromatic
  polynomial}, Algebr. Geom. Topol. \textbf{5} (2005), 1365--1388 (electronic).
  \MR{MR2171813 (2006g:57020)}

\bibitem[Kau87]{Kauffman}
Louis~H. Kauffman, \emph{State models and the {J}ones polynomial}, Topology
  \textbf{26} (1987), no.~3, 395--407. \MR{MR899057 (88f:57006)}

\bibitem[Kho00]{Khovanov}
Mikhail Khovanov, \emph{A categorification of the {J}ones polynomial}, Duke
  Math. J. \textbf{101} (2000), no.~3, 359--426. \MR{MR1740682 (2002j:57025)}

\bibitem[KR05]{Rozansky}
Mikhail Khovanov and Lev Rozansky, \emph{Matrix factorizations and link
  homology ii}, \textsf{arXiv:math/0505056}, 2005.

\bibitem[OS04]{Ozsvath}
Peter Ozsv{\'a}th and Zolt{\'a}n Szab{\'o}, \emph{Holomorphic disks and knot
  invariants}, Adv. Math. \textbf{186} (2004), no.~1, 58--116. \MR{MR2065507
  (2005e:57044)}

\bibitem[OS07]{Szabo}
\bysame, \emph{A cube of resolutions for knot floer homology},
  \textsf{arXiv:0705.3852}, 2007.

\bibitem[Vir02]{Viro}
Oleg Viro, \emph{Remarks on definition of khovanov homology},
  \textsf{arXiv:math.GT/0202199}, 2002.

\end{thebibliography}

\end{document}